\documentclass[12pt]{amsart}
\usepackage{amssymb,amsmath,amsfonts,latexsym}
\usepackage{bm}
\newcommand{\newsection}[1]{\setcounter{equation}{0} \section{#1}}
\def\A{{\underline{A}}}
\def\B{{\underline{B}}}

\newcommand{\bea}{\begin{eqnarray}}
\newcommand{\eea}{\end{eqnarray}}

\newcommand{\cld}{\mathcal{D}}
\newcommand{\cle}{\mathcal{E}}

\newcommand{\clg}{\mathcal{G}}
\newcommand{\clh}{\mathcal{H}}
\newcommand{\clk}{\mathcal{K}}
\newcommand{\cll}{\mathcal{L}}
\newcommand{\clm}{\mathcal{M}}
\newcommand{\cln}{\mathcal{N}}

\newcommand{\clp}{\mathcal{P}}
\newcommand{\clr}{\mathcal{R}}

\newcommand{\raro}{\rightarrow}

\def \qed {\hfill \vrule height6pt width 6pt depth 0pt}
\def\textmatrix#1&#2\\#3&#4\\{\bigl({#1 \atop #3}\ {#2 \atop #4}\bigr)}
\def\dispmatrix#1&#2\\#3&#4\\{\left({#1 \atop #3}\ {#2 \atop #4}\right)}
\newcommand{\be}{\begin{equation}}
\newcommand{\ee}{\end{equation}}
\newcommand{\ben}{\begin{eqnarray*}}
\newcommand{\een}{\end{eqnarray*}}

\newcommand{\NI}{\noindent}

\newcommand{\bi}{\begin{itemize}}
\newcommand{\ei}{\end{itemize}}

\def\5{{5\superprime}}

\newtheorem{Theorem}{\sc Theorem}[section]
\newtheorem{Lemma}[Theorem]{\sc Lemma}
\newtheorem{Proposition}[Theorem]{\sc Proposition}
\newtheorem{Corollary}[Theorem]{\sc Corollary}
\newtheorem{Definition}[Theorem]{\sc Definition}
\newtheorem{Example}[Theorem]{\sc Example}

\newtheorem{Note}[Theorem]{\sc Note}
\newtheorem{Question}{\sc Question}
\newtheorem{ass}[Theorem]{\sc Assumption}

\theoremstyle{remark}
\newtheorem{Remark}[Theorem]{\sc Remark}

\newcommand{\bt}{\begin{Theorem}}
\def\beginlem{\begin{Lemma}}
\def\beginprop{\begin{Proposition}}
\def\begincor{\begin{Corollary}}
\def\begindef{\begin{Definition}}
\def\beginexamp{\begin{Example}}
\def\beginrem{\begin{Remark}}
\def\beginq{\begin{Question}}
\def\beginass{\begin{ass}}
\def\beginnote{\begin{Note}}
\newcommand{\et}{\end{Theorem}}
\def\endlem{\end{Lemma}}
\def\endprop{\end{Proposition}}
\def\endcor{\end{Corollary}}
\def\enddef{\end{Definition}}
\def\endexamp{\end{Example}}
\def\endrem{\end{Remark}}
\def\endq{\end{Question}}
\def\endass{\end{ass}}
\def\endnote{\end{Note}}

\begin{document}

\title[Contractions with Polynomial characteristic functions]{Contractions with Polynomial characteristic functions I. Geometric approach}
\author[Foias]{Ciprian Foias}
\address[Ciprian Foias]{Department of Mathematics, Texas A\&M University, College Station, Texas 77843, USA}

\author[Sarkar]{Jaydeb Sarkar}
\address[Jaydeb Sarkar]{Department of Mathematics, The University of Texas at San Antonio, San Antonio, TX 78249, USA}
\email{jaydeb.sarkar@utsa.edu}

\thanks{This research was partially supported by a grant from the National Science Foundation. Part of this research was done while the second author was at Texas A\&M University, College Station, USA}


\begin{abstract}
In this note we study the completely non unitary contractions on separable complex Hilbert spaces which have polynomial characteristic functions. These operators are precisely those which admit a matrix representation of the form
\begin{equation*}
T = \begin{bmatrix}S & * & *\\0 & N & *\\0& 0& C \end{bmatrix},
\end{equation*}
where $S$ and $C^*$ are unilateral shifts of arbitrary multiplicities and $N$ is nilpotent. We prove that dimension of ker$S^*$ and dimension of $\mbox{ker}\,C$ are unitary invariants of $T$ and that $N$, up to a quasi-similarity is uniquely determined by $T$. Also, we give a complete classification of the subclass of those contractions for which their  characteristic functions are monomials.
\end{abstract}

\today

\maketitle
\newsection{Introduction}

One of the central problems in operator theory is the classification (up to unitary equivalence) of bounded linear operators on an infinite dimensional separable Hilbert space. In all this generality the problem may not have a satisfactory solution. Therefore, the focus of research on this problem was (and it still is) on classes of operators for which one can identify useful (unitary) invariants (e.g.  \cite{B}, \cite{CD}, \cite{GGK}, \cite{GK}, \cite{NF}, \cite{Pe} and  \cite{Pi}). The topic of this article is the class of all c.n.u. contractions for which the characteristic functions are operator-valued polynomials.

To be more specific, we recall the following basic concepts. A contraction $T$ on a Hilbert space $\clh$, (i.e. $\|T h\| \leq \|h\|$ for all $h$ in $\clh$) is completely non unitary (or also, c.n.u.) if there is no non-trivial $T$-reducing subspace $\clm$ of $\clh$ such that $T|_{\clm}$ is unitary. An operator valued analytic function $\Theta \colon\ \mathbb{D} \raro \cll(\cle, \cle_*)$ for some Hilbert spaces $\cle$ and $\cle_*$ is said to be contractive if $$\|\Theta(z) h \| \leq \|h\|, \quad \quad h \in \cle,$$ and purely contractive when it satisfies also
$$\|\Theta(0) f\| < \|f\|, \quad \quad f \in \cle \, \mbox{and}\, f \neq 0.$$

\NI It is known (see \cite{NF}, Chapter VI) that there is a bijective correspondence between the class of purely contractive operator valued analytic functions ``modulo coincidence'' and the class of c.n.u. contractions ``modulo unitary equivalence''. More precisely, for a given contraction $T$, one can define the defect operators $D_T = (I_{\clh} - T^* T)^{\frac{1}{2}}$ and $D_{T^*} = (I_{\clh} - T T^*)^{\frac{1}{2}}$ with defect spaces  $\cld_T = \overline{\mbox{Ran} D_T}$ and $\cld_{T^*} = \overline{\mbox{Ran} D_{T^*}}$. Then the characteristic function  of the contraction $T$ is the $\cll(\cld_T, \cld_{T^*})$-valued purely contractive analytic function defined by
\begin{equation}\label{ch-def}
\Theta_T (z) = [ - T + z D_{T^*} (I_{\clh} - z T^*)^{-1} D_{T}]|_{\cld_T},
\end{equation}
\NI for all $z$ in $\mathbb{D}$. Two operator valued analytic functions $\Theta \colon\  \mathbb{D} \raro \cll(\clm, \clm_*)$ and $\Psi : \mathbb{D} \raro \cll(\cln, \cln_*)$ are said to coincide if there exist two unitary operators $\tau \colon\ \clm \raro \cln$ and $\tau_* \colon\ \clm_* \raro \cln_*$ such that $\tau_* \Theta(z) = \Psi(z) \tau$. Two c.n.u. contractions $T$ on $\clh$ and $R$ on $\clk$ are unitarily equivalent (that is, there is a unitary operator $U$ from $\clh$ to $\clk$ such that $T = URU$), if and only if the characteristic
functions $\Theta_T(z)$ and $\Theta_R(z)$ coincide. The converse of this fact is also true (see \cite{NF}, Theorem 3.4). Moreover, for a given $\cll(\cle, \cle_*)$-valued purely contractive analytic function $\Theta(z)$ defined on $\mathbb{D}$, there exists a c.n.u. contraction $T$ on some Hilbert space such that $\Theta_T(z)$ coincide with $\Theta(z)$.

In this paper we prove that a c.n.u. contraction $T$ has a polynomial characteristic function of degree $n$ if and only if $T$ admits a matrix representation
\begin{equation}\label{1}
T = \begin{bmatrix}S & * & *\\0 & N & *\\0& 0& C \end{bmatrix},
\end{equation}
where $S$ and $C^*$ are unilateral shifts of arbitrary multiplicities and $N$ is a nilpotent operator of order $n$. We stress that the spaces on which one (or two) of the diagonal entries of (\ref{1}) acts may be $\{0\}$ and the representation (\ref{1}) becomes one of the following $$T = \begin{bmatrix} S& *\\ 0& N\end{bmatrix}, \begin{bmatrix} S& *\\ 0& C\end{bmatrix}, \begin{bmatrix} N& *\\ 0& C\end{bmatrix},\quad \mbox{or}\quad \begin{bmatrix} S\end{bmatrix}, \begin{bmatrix} N\end{bmatrix}, \begin{bmatrix} C\end{bmatrix}.$$ For convenience, we will refer to these representations as degenerate forms of (\ref{1}).

One of our aims is to identify some properties of the diagonal entries in (\ref{1}) which are unitary invariant of $T$.

In the next section, we study some relatively simple examples which inspired our research, namely  the c.n.u. contractions for which the  characteristic function is a scalar polynomial of degree $ \leq 1$. We point out that even for this simple case there are natural questions which are not trivial.

In Section 3, we give the results already mentioned before, namely - the c.n.u. contractions with a polynomial characteristic function are precisely those which have a upper triangular representation (\ref{1}). Moreover, we define and study two canonical such representations.

In Section 4 we show that the multiplicities of the shift $S$ and the co-shift $C$ in the triangular representation (\ref{1}) are unitary invariants of $T$.

In Section 5 we define the "minimal" versions of the canonical upper triangular representations of type (\ref{1}) considered in Section 3 and show that their "central" nilpotent entries are quasi-similar.

In Section 6 we study the degenerate forms of the representation (\ref{1}), and characterize them in term of the characteristic function of $T$.

Then, in Section 7, we consider the special case of the monomial characteristic functions.

\NI\textsf{Acknowledgement:} We thank the referee for his observations and suggestions. In particular, the general result in Step 2 of Theorem \ref{S-C} is due to the referee. Our original proof was less transparent.

\newsection{An Illustrative Example}

In this section, we give a concrete example for which the resolution of our general problem presented in the introduction can be obtained by elementary computations. In particular, this example serves as a paradigm for our general considerations, since it provides a complete picture of the c.n.u. contractions for which the characteristic functions coincides with a scalar polynomial of degree one.

Let $l^2(\mathbb{Z})$ be the Hilbert space of all square summable sequences given by
$$
l^2(\mathbb{Z}) = \left\{ (\alpha_n) = (\ldots, \alpha_2, \alpha_1, \alpha_0, \alpha_{-1}, \alpha_{-2}, \ldots) \colon \  \alpha_k \;\;\text{in}\;\; \mathbb{C} \;\; \text{and} \;\; \sum_{k \in \mathbb{Z}} |\alpha_k|^2 < \infty\right\},
$$
with the standard orthonormal basis $\{e_n, e_0, e_{-n}\}_{n \geq 1}$. First we construct a bounded linear operator on $l^2(\mathbb{Z})$ as follows:

\begin{Definition}\label{definitionTabc}
Let $a, b, c$ be in $\bar{\mathbb{D}}$ and $T_{a,b,c}$ be the bounded linear operator on $l^2(\mathbb{Z})$ defined by

\begin{equation}\label{T_abc}
T_{a,b,c} e_n =  \left\{
\begin{array}{cl}
a e_1&\quad\text{if }n = 0\\
b e_0 + c e_1 & \quad \text{if } n = -1\\
e_{n+1} & \quad \text{if } n \neq 0, -1.
\end{array}
\right.
\end{equation}
\end{Definition}
\vspace{0.1in}
The adjoint of $T_{a.b.c}$ is given by
\begin{equation}
T_{a,b,c}^* e_n =  \left\{
\begin{array}{cl}
\bar{a} e_0 + \bar{c} e_{-1}&\quad\text{if }n = 1\\
\bar{b} e_{-1}& \quad \text{if } n = 0\\
e_{n-1} & \quad \text{if } n \neq 1, 0.
\end{array}
\right.
\end{equation}

It is easy to see that with respect to the decomposition of $l^2(\mathbb{Z})$ as $$l^2(\mathbb{Z}) = \clh_1 \oplus \clh_0 \oplus \clh_{-1},$$ where $\clh_1 = \mbox{span}\{e_1, e_2, \ldots \}$, $\clh_0 = \mbox{span}\{e_0\}$ and $\clh_{-1} = \mbox{span} \{ e_{-1}, e_{-2}, \ldots\}$, $T_{a,b,c}$ can be expressed as the upper triangular matrix

\begin{equation}\label{2.3}
T_{a, b, c} = \left[
\begin{array}{ccc}
S & * & *\\
0 & N & *\\
0& 0& C  \end{array} \right] = \left[
  \begin{array}{ c c c c|c|c c c c}
\vdots&\vdots &\vdots&\vdots   &\vdots  &\vdots& \vdots&\vdots&\vdots\\
\vdots&\bm{0}&1&0  &0&  0&0&0&\vdots\\
\vdots&0&\bm{0}&1  &0&  0&0&0&\vdots\\
\vdots&0&0&\bm{0}  &a&  c&0&0&\vdots\\ \hline

\vdots&0&0&0   &\bm{0}&   b&0&0&\vdots\\ \hline

\vdots&0&0&0   &0&  \bm{0}& 1&0&\vdots\\
\vdots&0&0&0   &0&  0&\bm{0}&1&\vdots\\
\vdots&0&0&0  &0&   0& 0&\bm{0}&\vdots\\
\vdots&\vdots&\vdots&\vdots      &\vdots&   \vdots&\vdots&\vdots&\bm{\ddots}\\

\end{array} \right],
\end{equation}

\NI where $S$ is the shift of multiplicity one on $\clh_1$, $N$ ($= 0$) is the nilpotent of order one on $\clh_0$ and $C$ is the co-shift of multiplicity one on $\clh_{-1}$.

\NI Observe that $T_{a,b,c}$ is a contraction if and only if the matrix
\[
Q_{a,b,c} = \left[
\begin{array}{cc}
 a&c\\0&b
\end{array}
\right]
\]
is a contraction. It is elementary to check that $Q_{a,b,c}$ is a contraction if and only if
\begin{equation}\label{abcgamma1}
|a|^2 + |c|^2, \;\; |b|^2 + |c|^2 \leq 1
\end{equation}
and

\begin{equation}\label{abcgamma2}
c = \gamma (1-|a|^2)^{\frac{1}{2}} (1-|b|^2)^{\frac{1}{2}},
\end{equation}
where $\gamma$ is in $\bar{\mathbb{D}}$. We assume the conditions (\ref{abcgamma1}) and (\ref{abcgamma2}) on $a, b$ and $c$ so that the operator $T_{a,b,c}$ is a contraction.

Notice that $T_{a,b,c}$ is the bilateral shift when $|a|=|b|=1$ and $c =0$. Thus to avoid this special case, \textit{in what follows we will assume that} $\{a, b\} \cap \mathbb{D} \neq \phi$.

The squares of the defect operators of the contraction $T_{a,b,c}$ are given by
\begin{equation}\label{dt}
(I - T_{a,b,c}^* T_{a,b,c}) e_n =  \left\{
\begin{array}{cl}
0 & \text{if } n \neq 0, -1\\
(1 - |a|^2) e_0 - a\bar{c} e_{-1} & \mbox{if } n =0\\
- \bar{a}c e_0 + (1 - |b|^2 - |c|^2) e_{-1}  & \mbox{if } n = -1,
\end{array}
\right.
\end{equation}
and
\begin{equation}\label{dt*}
(I - T_{a,b,c} T_{a,b,c}^*) e_n =  \left\{
\begin{array}{cl}
0 & \text{if } n \neq 1, 0\\
(1 - |a|^2 - |c|^2) e_1 - b\bar{c} e_0 & \mbox{if } n = 1\\
- \bar{b}c e_1 + (1 - |b|^2) e_0 & \mbox{if } n = 0.
\end{array}
\right.
\end{equation}

\begin{Lemma}
The defect spaces $\cld_{T_{a,b,c}}$ and $\cld_{T_{a,b,c}^*}$ of $T_{a,b,c}$ are one dimensional if and only if
\begin{equation}\label{gamma1}
(1-|a|^2) (1-|b|^2)(1-|\gamma|^2) =0.
\end{equation}
In particular, if $a,b$ are in $\mathbb{D}$, then (\ref{gamma1}) is equivalent to $|\gamma| = 1$ in the representation (\ref{abcgamma2}) of $c$.
\end{Lemma}
\NI\textsf{Proof.} Due to (\ref{dt}) and (\ref{dt*}), the defect spaces of $T_{a, b, c}$ are of dimension one if and only if the determinants of the matrices

\[
\left[
\begin{array}{ccc}
1-|a|^2 & - \bar{a}c\\
-a \bar{c} & 1-|b|^2-|c|^2\end{array} \right] \quad \mbox{and} \quad \left[
\begin{array}{ccc}
1-|a|^2 -|c|^2 & - b\bar{c}\\
-\bar{b}c & 1-|b|^2\end{array} \right]
\]
are zero, which is equivalent to the condition (\ref{gamma1}). \qed
\vspace{0.1in}

According to this last result $T_{a,b,c}$ will have a scalar characteristic functions if and only if (\ref{gamma1}) holds. Therefore, \textit{throughout this Section 2 we will assume that (\ref{gamma1}) holds}.

 Now we study the special case when $c=0$. In this case, by (\ref{abcgamma2}) we have either $|a|=1$ or $|b|=1$. Let $|a| <1$ and $|b|=1$. By (\ref{dt}) and (\ref{dt*}) it follows that $\cld_{T_{a,b,0}} = \mathbb{C} e_0$ and $\cld_{T^*_{a,b,0}} = \mathbb{C} e_1$. Moreover, $\Theta_{T_{a,b,0}}(z) = a I_{\mathbb{C}e_0 \raro \mathbb{C}e_1}$ for all $z$ in $\mathbb{D}$, where $I_{\mathbb{C}e_0 \raro \mathbb{C}e_1}$ maps $e_0$ into $e_1$. Thus $\Theta_{T_{a,b,0}}(z)$ coincides with the scalar constant function $\Theta(z) = a$ for all $z$ in $\mathbb{D}$. Similarly, if $|a|=1$ then $\Theta_{T_{a,b,0}}(z)$ coincides with the scalar constant function $\Theta(z) = b$ for all $z$ in $\mathbb{D}$. We have thus obtained the following proposition.

\begin{Proposition}\label{ab0}
Let $\text{max}\{|a|,|b|\} = 1$ or equivalently, $c = 0$ and $\{|a|,|b|\} \cap \mathbb{D} \neq \phi$. Then the characteristic function of $T_{a,b,0}$ coincides with the scalar characteristic function $\Theta(z) = m  = \mbox{min}\{|a|, |b|\}$ for all $z$ in $\mathbb{D}$.
\end{Proposition}

Thus we have the following complement to Proposition \ref{ab0}.
\begin{Remark}
\textnormal{(i) In case $m >0$, $(\frac{a}{|a|}, \frac{b}{|b|})$ can be any fixed element of $\mathbb{T}^2$ and all $T_{a,b,0}$ are unitarily equivalent with $T_{|a|,|b|,0}$.\\
\;\;(ii) In case $m=0$, then if  $a \neq0$, $T_{a,0,0}$ is unitarily equivalent to $T_{|a|,0,0}$ and $\frac{a}{|a|}$ can be any fixed element of $\mathbb{T}$ while if $b \neq 0$, the similar fact holds for $T_{0,b,0}$.}
\end{Remark}

We shall now consider the case when $c \neq 0$.
\begin{Proposition}\label{thetaabc}
Let $a,b$ and $c$ satisfy the condition (\ref{abcgamma2}) with $|\gamma|=1$ and $c \neq 0$. Then $$\Theta_{T_{a,b,c}}(z) = (-ab) + c z,$$ or equivalently, $$\Theta_{T_{a,b,c}}(z) = (-ab) + \gamma (1-|a|^2)^{\frac{1}{2}} (1-|b|^2)^{\frac{1}{2}}z.$$
\end{Proposition}
\NI \textsf{Proof.} We conclude from (\ref{abcgamma2}), (\ref{dt}) and (\ref{dt*}) that

\begin{equation*}
D^2_{T_{a,b,c}} e_n =  \left\{
\begin{array}{cl}
0 & \text{if } n \neq 0, -1\\
\frac{ -(1 - |a|^2)}{c} (a(1 - |b|^2) e_{-1} - c e_0) & \mbox{if } n = 0\\
\bar{a} (a(1 - |b|^2) e_{-1} - c e_0) & \mbox{if } n =-1
\end{array}
\right.
\end{equation*}
and
\begin{equation*}
D^2_{T_{a,b,c}^*} e_n =  \left\{
\begin{array}{cl}
0 & \text{if } n \neq 1, 0\\
\frac{-b \bar{c}}{ 1 - |b|^2} ( (1 - |b|^2) e_0 - \bar{b} c e_1) & \mbox{if } n = 1\\
(1 - |b|^2) e_0 - \bar{b} c e_1& \mbox{if } n =0.
\end{array}
\right.
\end{equation*}

In particular, the defect spaces are given by $\cld_{T_{a,b,c}} = \mathbb{C} u$ and $\cld_{T_{a,b,c}^*} = \mathbb{C}v$, where $u = a(1 - |b|^2) e_{-1} - c e_0$ and $v = (1 - |b|^2) e_0 - \bar{b} c e_1$. Moreover, we have
$$T_{a,b,c} u = (- a b) v, D_{T_{a,b,c}} u = \sqrt{1 - |a|^2 |b|^2} u, D_{T_{a,b,c}^*} v = \sqrt{1 - |a|^2 |b|^2} v,$$
and $$D_{T_{a,b,c}^*} u = \frac{c}{\sqrt{1 - |a|^2 |b|^2}} v.$$

Thus $D_{T_{a,b,c}^*} D_{T_{a,b,c}} u = c v$ and hence by virtue of (\ref{ch-def}), the characteristic function  of $T_{a,b,c}$ is given by $$\Theta_{T_{a,b,c}}(z) = (-ab) + c z.$$\qed

\vspace{0.1in}
It is elementary to check the following characterization of a purely contractive scalar polynomial of degree one.
\begin{Lemma}\label{thetaaplha}
Let $\alpha, \beta \in \mathbb{C}$. Then the one degree polynomial $$\Theta(z) = \alpha + \beta z \qquad \quad (z \; \mbox{in}\; \mathbb{D})$$ is purely contractive if and only if
\begin{equation}\label{alphabeta}|\alpha| + |\beta| \leq 1 \quad \quad \text{and} \quad \quad |\alpha| < 1.\end{equation}
\end{Lemma}

\vspace{0.1in}

\begin{Remark}
\textnormal{Since each characteristic function is pure, so $\alpha$ and $\beta$ given by $\alpha = - a b$ and $\beta =c$ with $a,b,c$ as in Proposition \ref{ab0} or Proposition \ref{thetaabc} will satisfy the conditions in (\ref{alphabeta}). Conversely if two complex numbers $\alpha$ and  $\beta$ satisfy (\ref{alphabeta}) then it is easy to check that there exist two real numbers $a$ and $b$ satisfying $\alpha = -ab$ and (\ref{abcgamma2}) with $c=\beta$ and $|\gamma|=1$.}
\end{Remark}

Since the unitary part of $T_{a,b,c}$ can not be connected to the characteristic function  $\Theta_{T_{a,b,c}}(z)$, we need to determine when that unitary part is not present, that is, when $T_{a,b,c}$ is c.n.u. The answer to this question is given by the following.
\begin{Proposition}\label{abccnu}
Let $a,b$ and $c$ satisfy the condition (\ref{abcgamma2}) with $|\gamma| = 1$ and $c \neq0$. Then the contraction $T_{a,b,c}$ is c.n.u. if and only if $|c| <1$.
\end{Proposition}

\NI\textsf{Proof.} If $|c| = 1$ then due to (\ref{abcgamma1}) we have $a=b=0$. So $T_{a,b,c}$ is the direct sum of the operator $0$ on $\clh_0$ and the bilateral shift $U$ defined on $\clh_1 \oplus \clh_{-1}$ by
\begin{equation*}
Ue_n =  \left\{
\begin{array}{cc}
e_{n+1}&\;\;\text{if }n \neq -1\\
c e_1& \quad \text{if } n = -1.
\end{array}
\right.
\end{equation*}
Thus the contraction $T_{a,b,c}$ is not c.n.u.

Conversely, assume $|c|<1$. Consider a $T_{a,b,c}$-reducing subspace $\clm$ of $l^2(\mathbb{Z})$ such that $T_{a,b,c}|_{\clm}$ is unitary operator and let $(\alpha_k)$ be a vector in $\clm$.

\textsf{Case 1} (when $a=0$): We know from (\ref{T_abc}) that $\alpha_0=0$. As $T^{*l}_{a,b,c} (\alpha_k)$ is in $\clm$ for all $l \geq 1$, (\ref{dt}) says precisely that $\alpha_{l} =0$ for all $l \geq1$. Also $\alpha_l =0$ for all $l < 0$ follows from the fact that $T^{l}_{a,b,c} (\alpha_k)$ is in $\clm$ for all $l \geq 0$. Hence $\clm = \{\bm{0}\}$.

\textsf{Case 2} (when $b=0$): This follows from Case 1 by considering $T^*_{a,b,c}$ instead of $T_{a,b,c}$.

\textsf{Case 3} (when all $a,b,c \neq 0$): By (\ref{dt*}) and (\ref{dt}) we have
$$
\alpha_1 = \left(\frac{1-|b|^2}{1-|a|^2}\right)^{\frac{1}{2}} \frac{\gamma}{b} \alpha_0 \quad \quad \text{and} \quad \quad \alpha_{-1} = \left(\frac{1-|a|^2}{1-|b|^2}\right)^{\frac{1}{2}} \frac{\bar{\gamma}}{\bar{a}} \alpha_0.
$$

By repeated application of $T^*_{a,b,c}$ and $T_{a,b,c}$ on $(\alpha_k)$ we obtain that
$$
\alpha_{k+1} = \left(\frac{1-|b|^2}{1-|a|^2}\right)^{\frac{k+1}{2}} \left(\frac{\gamma \bar{a}}{b}\right)^{k+1} \frac{\alpha_0}{\bar{a}} \quad \quad \text{and}\quad \quad \alpha_{-(k+1)} = \left(\frac{1-|a|^2}{1-|b|^2}\right)^{\frac{k+1}{2}} \left(\frac{\bar{\gamma} b}{a}\right)^{k+1} \frac{\alpha_0}{b},
$$
for all $k \geq 1$. Therefore for all $k \geq1$
$$
|\alpha_{k+1}| = \left(\frac{1-|b|^2}{1-|a|^2}\right)^{\frac{k+1}{2}} |\frac{a}{b}|^{k+1} |\frac{\alpha_0}{a}| \quad \quad \text{and}\quad \quad |\alpha_{-(k+1)}| = \left(\frac{1-|a|^2}{1-|b|^2}\right)^{\frac{k+1}{2}} |\frac{b}{a}|^{k+1} |\frac{\alpha_0}{b}|.
$$
This with Cauchy-Schwarz inequality gives
$$
N \frac{|\alpha_0|^2}{|a| |b|} = \sum_{k=2}^{N} |\alpha_k| |\alpha_{-k}| \leq \left(\sum_{k=2}^{N} |\alpha_k|^2\right)^{\frac{1}{2}} \left(\sum_{k=2}^{N} |\alpha_{-k}|^2\right)^{\frac{1}{2}},
$$
for all $N \geq 2$ and hence since $(\alpha_k)$ is in $l^2(\mathbb{Z})$, we must have $\alpha_0 =0$ and so $(\alpha_k) = \bm{0}$. Therefore, $\clm = \{\bm{0}\}$ and this completes the proof. \qed

\vspace{0.1in}

\begin{Remark}\label{remabc}
\textnormal{Proposition \ref{thetaabc} and \ref{abccnu} entail that if the polynomial $\Theta(z)$ in Lemma \ref{thetaaplha} is a characteristic function of a c.n.u. contraction $T_{a, b, c}$ then $\beta \in \mathbb{D}$.}
\end{Remark}

Summarizing the previous results we obtain the following.

\begin{Theorem}\label{scalar}
The scalar polynomial $\Theta(z) = \alpha + \beta z $ is a characteristic function of a c.n.u. contraction if and only if $|\alpha| <1$ and $|\alpha| + |\beta| \leq1$. In this case a c.n.u. contraction for which the characteristic function coincides with $\Theta$ is unitarily equivalent to any contraction $T_{a, b, c}$ defined by the equation (\ref{T_abc}), where $c = \beta$ and $a, b \in \overline{\mathbb{D}}$ satisfy  $ab = - \alpha $ and $(1-|a|^2) (1-|b|^2) = |\beta|^2$.
\end{Theorem}
\vspace{0.2in}

A natural question related to Theorem \ref{scalar} is whether the diagonal terms in a triangular representation of
\begin{equation*}
T = T_{a,b,c} = \left[
\begin{array}{ccc}
S & * & *\\
0 & N & *\\
0& 0& C  \end{array} \right],
\end{equation*}
where $S$ and $C^*$ are shifts and $N$ is a nilpotent operator of order one are unitary invariants of $T$. To have a good indication that the question is not trivial, let us consider the case $\alpha = 0$, or to be more specific, the case $a = 0$. Then $b \neq 0$ and the operator $T_{a,b,c}$ has the triangular representation (\ref{2.3}) and as well the following

\begin{equation}
T_{a, b, c} = \left[
\begin{array}{ccc}
S & * & *\\
0 & N & *\\
0& 0& C  \end{array} \right] = \left[
  \begin{array}{ c c c c|cc|c c c c}
\vdots&\vdots &\vdots&\vdots   &\vdots&\vdots  &\vdots& \vdots&\vdots&\vdots\\
\vdots&\bm{0}&1&0  &0&0&  0&0&0&\vdots\\
\vdots&0&\bm{0}&1  &0&0&  0&0&0&\vdots\\
\vdots&0&0&\bm{0}  &1&0&  0&0&0&\vdots\\ \hline
\vdots&0&0&\bm{0}  &\bm{0}&0&  c&0&0&\vdots\\

\vdots&0&0&0   &0&\bm{0}&   b&0&0&\vdots\\ \hline
\vdots&0&0&0   &0&0&   \bm{0}&1&0&\vdots\\

\vdots&0&0&0   &0&0&  0&\bm{0}&1&\vdots\\
\vdots&0&0&0   &0&0&  0&0&\bm{0}&\vdots\\
\vdots&\vdots&\vdots&\vdots      &\vdots&\vdots&   \vdots&\vdots&\vdots&\bm{\ddots}\\

\end{array} \right],
\end{equation}
for which the shift and the co-shift are unitarily equivalent to their analogues in (\ref{2.3}) but the two "central" nilpotent operators are
\[
\left[ \begin{array}{c}
0\end{array}\right] \quad \quad \mbox{and} \quad \quad \left[ \begin{array}{cc}
0 & 0\\
0& 0\end{array}\right],
\]
which are not even similar.

In what follows, we will study the general variant of the above question for the class of c.n.u. contractions with operator valued polynomial characteristic functions.


\newsection{Polynomial characteristic functions}

To start the study of the c.n.u. contractions with polynomial characteristic functions, we need to  recall that the order of a nilpotent operator $N$ is the smallest power $p$ for which $N^p = 0$.

\begin{Proposition}\label{SNC}
Let $\clh$ be a Hilbert space with the decomposition $\clh = \clh_1 \oplus \clh_0 \oplus \clh_{-1}$. Also let $S$ in $\cll(\clh_1)$ be an isometry, $N$ in $\cll(\clh_0)$ be a nilpotent of order $n$ and $C$ be a co-isometry in $\cll(\clh_{-1})$. Then the characteristic function  of any contraction of the form
\begin{equation*}
T = \begin{bmatrix}S & * & *\\0 & N & *\\0& 0& C \end{bmatrix}
\end{equation*}
is a polynomial of degree $ \leq n$.
\end{Proposition}
\NI\textsf{Proof.} It is easy to see that
\begin{equation*}
D^2_T = \begin{bmatrix}
0 & * & * \\ * & * & * \\ *& *& *
\end{bmatrix}.
\end{equation*}
Let
\begin{equation*}
D_T = \begin{bmatrix}
A & B & C\\B^* & * & *\\C^*& *& *
\end{bmatrix},
\end{equation*}
where $A$ is a self adjoint operator. Then,

\begin{equation*}
D^2_T = \begin{bmatrix}
A^2+BB^*+CC^* & * & *\\ * & * & *\\ *& *& *
\end{bmatrix} = \begin{bmatrix}
0 & * & *\\ * & * & *\\ *& *& *
\end{bmatrix}.
\end{equation*}
Thus the sum of positive operators $A^2+BB^*+CC^*$ is the zero operator and consequently $A=0$, $B=0$ and $C=0$, and
\begin{equation*}
D_T = \begin{bmatrix}
0 & 0 & 0\\0 & * & *\\0& *& *
\end{bmatrix}.
\end{equation*}
In a similar way we deduce
\begin{equation*}
D_{T^*} = \begin{bmatrix}
* & * & 0\\ * & * & 0\\0& 0& 0
\end{bmatrix}.
\end{equation*}
Moreover,
\begin{equation*}
T^{*k} = \begin{bmatrix}
S^{*k} & 0 & 0\\ * & N^{*k} & 0\\ *& *& C^{*k}\end{bmatrix}
\end{equation*}

\NI so that $D_{T^*} T^{*k} D_T = 0$ for all $k \geq n$. Thus from \eqref{ch-def} we obtain
\[
\Theta_T(\lambda) = [-T + \sum^n_{k=1} \lambda^kD_{T^*} T^{*k-1}D_T] \ \Big| \ \cld_T
\]
which is a polynomial of degree $\le n$.\qed

\vspace{0.3in}

\begin{Remark}
If in Proposition \ref{SNC} we consider the representation of the operator $T_{a,b,c}$ with ${a, b, c} \in \mathbb{D} \setminus \{0\}$ in which $\clh_0 = \mathbb{C} e_1 + \mathbb{C} e_0 + \mathbb{C} e_{-1}$ then $N$ is of order 3. Although the characteristic function is a polynomial of degree one. Thus in Proposition \ref{SNC}, the degree of the characteristic function can be less than the order of $N$.
\end{Remark}

Next we consider a converse of Proposition \ref{SNC}.

\begin{Theorem}\label{prop1}
Let $T$ be a c.n.u. contraction on $\clh$ such that the characteristic function $\Theta_T$ is a polynomial of degree $n$. Then there exists three subspaces $\clh_1, \clh_0, \clh_{-1}$ of $\clh$ such that $\clh = \clh_1 \oplus \clh_0 \oplus \clh_{-1}$ and with respect to that decomposition, $T$ admits the matrix representation
\begin{equation}\label{3.1}
T = \begin{bmatrix}S & * & *\\0 & N & *\\0& 0& C
\end{bmatrix}
\end{equation}
where $S$ in $\cll(\clh_1)$ is a pure isometry, $N$ in $\cll(\clh_0)$ is a nilpotent of order $n$ and $C$ in $\cll(\clh_{-1})$ is a pure co-isometry; in particular, $S$ and $C^*$ are shifts of some multiplicities on $\clh_1$ and $\clh_{-1}$ respectively.

Furthermore in the class of representations \eqref{3.1} with the above specified properties there exist unique representations
\begin{alignat}{2}\label{eq3.2}
T &= \begin{bmatrix} S^{(c)}&*&*\\ 0&N^{(c)}&*\\ 0&0&C^{(c)}\end{bmatrix}\quad &\text{on} \quad \mathcal{H} &= \mathcal{H}^{(c)}_1 \oplus \clh^{(c)}_0 \oplus \mathcal{H}^{(c)}_{-1},\\
\label{eq3.3}
T &= \begin{bmatrix} S^{(c)}_*&*&*\\ 0&N^{(c)}_*&*\\ 0&0&C^{(c)}_*\end{bmatrix}\quad &\text{on}\quad \mathcal{H} &= \mathcal{H}^{(c)}_{1*} \oplus \mathcal{H}^{(c)}_{0*} \oplus \mathcal{H}^{(c)}_{-1*}
\end{alignat}
such that
\begin{align}\label{eq3.4}
 &\mathcal{H}^{(c)}_1 \subset \mathcal{H}_1,\quad \mathcal{H}^{(c)}_{-1} \supset \mathcal{H}_{-1},\\
\label{eq3.5}
&\mathcal{H}^{(c)}_{1*} \supset \mathcal{H}_1,\quad \mathcal{H}^{(c)}_{-1*} \subset \mathcal{H}_{-1}
\end{align}
holds for all other representations \eqref{3.1} in the class; moreover $\mathcal{H}^{(c)}_{-1}$, respectively $\mathcal{H}^{(c)}_{1*}$, is the maximal invariant subspace of $T^*$, respectively $T$, on which $T^*$, respectively $T$, is isometric.
\end{Theorem}
\vspace{0.1in}
\NI \textsf{Proof.} By the given assumption, it follows that $D_{T^*} T^{*(n+k)} D_T = 0$, which is equivalent to $D_{T} T^{(n+k)} D_{T^*} = 0$ for all $k \geq 0$. Define
$$
\clh^{(c)}_1 = \overline{\mbox{span}} \{ T^k h \colon\ h \in \mbox{in} \; \cld_{T^*}, k \geq n\}.
$$
It is easy to see that the subspace $\clh^{(c)}_1$ of $\clh$ is an invariant subspace of $T$ and $T|_{\clh_1}$ is a pure isometry. The last claim follows directly from the fact $T$ is c.n.u.\ and that $\|Th\| = \|h\|$ if and only if $D_T h = 0$ for $h$ in $\clh$. Now we construct another $T$ invariant subspace $\clm \supseteq \clh^{(c)}_1$  defined by
$$
\clm = \overline{\mbox{span}} \{ T^k h \colon\ h \in \cld_{T^*}, k \geq 0\},
$$
and a semi-invariant subspace of $T$ as $\clh^{(c)}_0 = \clm \ominus \clh^{(c)}_1$. We claim that the operator $N^{(c)} = P_{\clh^{(c)}_0} T|_{\clh^{(c)}_0}$ on $\clh^{(c)}_0$ is a nilpotent of order $n$. The claim follows once we observe that $N^{(c)k} = P_{\clh^{(c)}_0} T^k|_{\clh^{(c)}_0}$ for any $k \geq 0$,  that $T^{n} \clm \subseteq \clh^{(c)}_1$ and then by Theorem~\ref{SNC}. Therefore, we obtain that
\begin{equation*}
T|_{\clh^{(c)}_1 \oplus \clh^{(c)}_0} =
\begin{bmatrix}
S^{(c)} & * \\0 & N^{(c)}
\end{bmatrix},
\end{equation*}
where $S^{(c)} = T|_{\clh^{(c)}_1}$ is a pure isometry. Now we prove  that the operator $C^{(c)}$ on the $T$ co-invariant subspace $\clh^{(c)}_{-1} = \clh \ominus \clm$ defined by $C^{(c)} = P_{\clh^{(c)}_{-1}} T|_{\clh^{(c)}_{-1}}$  is a co-isometry.
Note that $h$ is in $\clh^{(c)}_{-1}$ if and only if $D_{T^*} T^{*k} h \perp \clh$ for all $k \geq 0$, or in other wards
\begin{equation}\label{eq3.6}
\clh^{(c)}_{-1} = \{ h \in \clh \colon\  D_{T^*} T^{*k} h = 0, k \geq 0\}.
\end{equation}
As a result, for $h$ in $\clh^{(c)}_{-1}$ and for $k \geq 0$, we have $\|T^{* (k+1)}\| = \|T^{*k} h\|$, or equivalently
\begin{equation}\label{clm_-1}
\clh^{(c)}_{-1} = \clm_{-1} := \{ h \in \clh \colon\ \|T^{*k+1} h\| = \| T^{*k} h\|, k \geq 0\},
\end{equation}
clearly the largest co-invariant subspace of $\clh$ such that the restriction of $T^*$ on it is isometric. It remains to prove the minimality property of $\clh^{(c)}_1$. For this let
\begin{equation}\label{rep*}
T = \begin{bmatrix}{S} & * & *\\ 0 & {N} & *\\ 0&0 & {C}
\end{bmatrix},
\end{equation}
be a representation of $T$ on $\clh = \clh'_1 \oplus \clh'_0 \oplus \clh'_{-1}$ such that $S'$ on $\clh'_1$ is an isometry, $N'$ on $\clh'_0$ is a nilpotent operator of order $n$ and $C'$ on $\clh'_{-1}$ is a co-isometry. Then the computation in the proof of Proposition \ref{SNC} yields
\begin{equation*}
D_{T^*} = \begin{bmatrix}* & * & 0\\ * & * & 0\\ 0&0 & 0
\end{bmatrix}.
\end{equation*}
Now observe that
\begin{equation*}
T^n D_{T^*} = \begin{bmatrix}{S^n} & * & *\\ 0 & {0} & *\\ 0&0 & {C^n}
\end{bmatrix} \begin{bmatrix}* & * & 0\\ * & {*} &0\\ 0&0 & 0 \end{bmatrix} = \begin{bmatrix}*& * & 0\\ 0 & {0} & 0\\ 0&0 &0
\end{bmatrix}.
\end{equation*}
Then we must have $T^{k+n} D_{T^*} \clh \subseteq \clh_1$ for all $k \geq 0$ and hence $\clh^{(c)}_1 \subseteq \clh_1$. This completes the proof concerning the representation \eqref{eq3.2}. For the representation \eqref{eq3.3}, the proof is similar once we define
\begin{equation}\label{eq3.9}
 \left\{\begin{array}{l}
\clh^{(c)}_{1*} = \{h\in \clh\colon \ \|T^{k+1}h\| = \|T^kh\|, k\ge 0\},\\
\clh^{(c)}_{-1*} = \overline{\text{span}}\{T^{*k}h\colon \ h\in \cld_T, k\ge m\},\\
\clh^{(c)}_{0*} = \clh \ominus (\clh^{(c)}_{1*} \oplus \clh^{(c)}_{-1*}).
        \end{array}
\right. \qed
\end{equation}

\vspace{0.2in}

Combining Proposition \ref{SNC} and Theorem \ref{prop1} we can readily obtain the following.
\begin{Corollary}\label{poly-class}
Let $T$ be a c.n.u. contraction on a Hilbert space $\clh$. Then the characteristic function $\Theta_T$ of $T$ is a polynomial if and only if there exists three closed subspaces $\clh_1, \clh_0, \clh_{-1}$ of $\clh$ with $\clh = \clh_1 \oplus \clh_0 \oplus \clh_{-1}$ and a pure isometry $S$ in $\cll(\clh_1)$, a nilpotent $N$ in $\cll(\clh_0)$ and a pure co-isometry $C$ in $\cll(\clh_{-1})$ such that $T$ has the following matrix representation
\begin{equation*}
T = \begin{bmatrix}
S & * & *\\0 & N & *\\0& 0& C
\end{bmatrix}.
\end{equation*}
Moreover the degree of $\Theta_T$ is the smallest order of $N$ possible in the matrix representations \eqref{eq3.4} of $T$.
\end{Corollary}

\vspace{0.2in}

\begin{Remark}\label{rem3.5}
Due to their uniqueness, it is convenient to call the representations \eqref{eq3.2} and \eqref{eq3.3} the \emph{canonical} and respectively  $*$-\emph{canonical} representations of the operator $T$ considered in Theorem~\ref{prop1}. Note that the latter can be obtained also from the canonical representation of $T^*$ by passing to the adjoint and by flipping the extreme spaces. In the sequel we will omit the superscripts $(c)$ (indicating that we are dealing with a canonical or $*$-canonical matrix representation) whenever no confusion can occur.
\end{Remark}

\begin{Remark}\label{rem3.6}
Let $T\in \cll(\clh), T'\in \cll(\clh')$ be c.n.u.\ contractions and assume that their characteristic functions are polynomials of degree $n$. Let
\begin{alignat*}{2}
 T &= \begin{bmatrix} S&*&*\\ 0&N&*\\ 0&0&C\end{bmatrix} \quad &\text{on}\quad \clh &= \clh_1 \oplus \clh_0 \oplus \clh_{-1},\\
T' &= \begin{bmatrix} S'&*&*\\ 0&N'&*\\ 0&0&C'\end{bmatrix} \quad &\text{on}\quad \clh' &= \clh'_1 \oplus \clh'_0 \oplus  \clh'_{-1}
\end{alignat*}
be their canonical representations. If $U\in \cll(\clh,\clh')$ is a unitary operator such that $UT=T'U$, then from the definitions \eqref{eq3.6}, \eqref{clm_-1} and \eqref{eq3.9} it readily follows that
\[
 U\clh_i = \clh'_i,\qquad (i=-1,0,1)
\]
and consequently
\[
 S'(U|\clh_1) = (U|\clh_1)S,\quad N'(U|\clh_0) = (U|\clh_0)N\quad \text{and}\quad C'(U|\clh_{-1}) = (U|\clh_{-1})C,
\]
that is, the diagonal entries $S$ and $S',N$ and $N',C$ and $C'$ are unitarily equivalent too. Obviously the similar fact holds true for $*$-canonical representations.
\end{Remark}

\begin{Remark}\label{rem3.7}
We point out that the canonical and the $*$-canonical representations of the contraction $T$ in Theorem~\ref{prop1} may be different, since a priori $\clh^{(c)}_1$ (resp.\ $\clh^{(c)}_{-1*}$) may not equal $\clh^{(c)}_{1*}$ (resp.\ $\clh^{(c)}_{-1}$). This is well illustrated by the properties of the operator $T_{a,b,c}$, considered in Section 2 (see also Theorem \ref{scalar}) presented below.
\end{Remark}

\begin{Proposition}
Let $a, b, c \in \mathbb{D}  \setminus \{0\}$ satisfy $c = \gamma (1 - |a|^2)^{\frac{1}{2}} ( 1 - |b|^2)^{\frac{1}{2}}$ with $|\gamma| = 1$ and let $T_{a,b,c}$ be the operator introduced in Definition \ref{definitionTabc}. Then the matrix representation (\ref{2.3}) of $T_{a,b,c}$ is canonical, respectively $*$-canonical if and only if $|a| \geq  |b|$, respectively $|b| \leq |a|$. In particular, the canonical and the $*$-canonical representations of $T_{a,b,c}$ coincide if and only if $|a| = |b|$ and in this case they coincide with (\ref{2.3}).
\end{Proposition}

\NI\textsf{Proof.} Let $(\alpha_k)$ be a vector in $\clh^{(c)}_{-1}  = \cap_{l=0}^{\infty} \mbox{ker~} D_{T^*_{a,b,c}} T^{* l}_{a,b,c} \supseteq \clh_{-1}$. As $D_{T_{a,b,c}}^* (\alpha_k) = 0$, by (\ref{dt*}) we have \[\alpha_1 = \left({\frac{1-|b|^2}{1-|a|^2}}\right)^{\frac{1}{2}} \frac{1}{b \bar{\gamma}} \alpha_0.\]By repeated application of $T^*_{a,b,c}$ on $(\alpha_k)$ and then using $D_{T^*_{a,b,c}} T^{*l}(\alpha_k) = 0$ we obtain that \[\alpha_k = \left[\left(\frac{1-|b|^2}{1-|a|^2}\right)^{\frac{1}{2}} \frac{\bar{a}}{b \bar{\gamma}}\right]^{k-1} \alpha_1,\]for all $k \geq 1$. As $\sum_{k=0}^{\infty} |\alpha_k|^2 < \infty$, we deduce that $\clh^{(c)}_{-1} \supsetneqq \clh_{-1}$ if and only if \[ |\left(\frac{1-|b|^2}{1-|a|^2}\right)^{\frac{1}{2}} \frac{\bar{a}}{b \bar{\gamma}}| < 1,\]or, equivalently \[|a| < |b|.\] Consequently,  $\clh^{(c)}_{-1} = \clh_{-1}$ if and only if $|a| = |b|$. Proceeding in a similar way with $\clh^{(c)}_{1*}$, one can obtain the proof of the remaining part of the statement.\qed


\newsection{Uniqueness of the shift and the co-shift}

In this section we discuss the uniqueness problem for the diagonal of the upper triangular representations of a c.n.u. contraction $T$ on $\clh$ with a polynomial characteristic function of degree $n$. It will be convenient to denote the family of these contractions by $\clp_n$.



First, we need to recall the following fact observed by Sz.-Nagy and the first author (see Remark 3.3, page 134 in \cite{FF}) and Parrott (see \cite{Pa}) on the isometric commutant lifting theorem.

Let $T$ in $\cll(\clh, \clk)$ and $U$ in $\cll(\clh', \clk')$ be two operators. Then we say that $U$ is a lifting of $T$ if $\clh' \supseteq \clh, \clk'\supseteq \clk$ and $P_{\clk} U = T P_{\clh}$, or equivalently, $U$ has a matrix representation of the form
\begin{equation}
U = \begin{bmatrix} T&0\\A&B\end{bmatrix} \colon\ \clh' = \clh \oplus \clh^{\perp} \raro \clk' = \clk \oplus \clk^{\perp}.
\end{equation}
In addition, if $U$ is an isometric operator (and hence $T$ is a contraction), then we say that $U$ is an isometric lifting of $T$. The minimal isometric dilation of a contraction is also an isometric lifting of the contraction. With this preliminary, the fact alluded above is the following. If $U$ on $\clk$ is an isometric lifting of a contraction $T$ in $\clh$, then $U$ admits a decomposition of the form $U = U_{m} \oplus U_{r}$ on $\clk = (\bigvee_{k=0}^{\infty} U^k \clh) \oplus \clk_r$ where $U_m$ is the minimal isometric dilation of the contraction $T$.

\begin{Theorem}\label{S-C}
Let $T$ be a c.n.u. contraction on $\clh$ such that the characteristic function $\Theta_T$ is a polynomial of degree $n$ and
\begin{equation*}
\begin{bmatrix} S&*&*\\0&N&*\\0&0&C\end{bmatrix}\quad \mbox{and} \quad \begin{bmatrix} {S_*}&*&*\\0&{N_*}&*\\0&0&{C_*}\end{bmatrix}
\end{equation*}
be the canonical and $*$-canonical matrix representations of $T$ on $\clh = \clh_1 \oplus \clh_0 \oplus \clh_{-1}$ and on $\clh = \clh_{1*} \oplus \clh_{0*} \oplus \clh_{-1*}$, respectively. Moreover, let
\begin{equation*}
T = \begin{bmatrix} S'&*&*\\0&N'&*\\0&0&C'\end{bmatrix},
\end{equation*}
be another representation of $T$ on $\clh = \clh'_{1} \oplus \clh'_{0} \oplus \clh'_{-1}$ where $S'$ is a pure shift on $\clh'_1$, $N'$ is a nilpotent operator of order $n$ on $\clh'_0$ and $C'$ is a pure co-shift on $\clh'_{-1}$. Then the pure shifts $S, S'$ and ${S_*}$ are unitarily equivalent and the pure co-shifts $C, C'$ and $C_*$ are unitarily equivalent.
\end{Theorem}
\NI \textsf{Proof.} It is enough to prove that $S_*$ and $S'$ are unitarily equivalent. The proof will be separated in three steps in which the result in Step 2 has an independent interest.

\NI\textsf{Step 1:} First observe that if we view $T$ as
\[T =  \begin{bmatrix}S'&A\\0&B \end{bmatrix}\]on $\clh'_1 \oplus \clh'_2$ where \[B = \begin{bmatrix}N'&*\\0&C'\end{bmatrix}\] on $\clh'_2 = \clh'_0 \oplus \clh'_{-1}$ then \[B^k \raro 0\]strongly as $k \raro \infty$. Indeed if \[C'^l h'_{-1} = 0,\]for some $l \,( = 1, 2, \ldots)$ and $h' = h'_0 \oplus h'_{-1}$, then we obviously have that \[B^k h'= 0,\]provided $k > n+l+1$. Thus \[\|B^k h'\| \raro o \quad \quad \mbox{for~} k\raro \infty\]
on a dense set in $\clh'_2$. But $\|B\| \leq 1$ and thus the claim follows.

\NI\textsf{Step 2:} Let $T \in \cll(\clh_1 \oplus \clh_2)$ be a c.n.u.  contractions such that \[T = \begin{bmatrix}S&A\\0&B\end{bmatrix}\]on $\clh_1 \oplus \clh_2$ where $S$ is an isometry and $B \in C_{0 \cdot}$ (that is, $B^k \raro 0$ strongly as $k \raro \infty$). Let moreover $\clm$ be the largest invariant subspace of $T$ on which the restriction of $T$ is isometric. Then
\begin{equation}\label{e4.2}\mbox{dim} \,\mbox{ker~}S^* = \mbox{dim~}\mbox{ker~} (T|_{\clm})^*.\end{equation}
To establish this fact let $\clg = \clm \ominus \clh_1$ and $X = P_{\clg} T|_{\clg}$. Then $T|_{\clm}$ is a c.n.u. isometric lifting of $X$ and \[X^* = (T|_{\clm})^*|_{\clg}.\] Thus $X^{*k} \raro 0$ strongly as $k \raro \infty$, that is, $X \in C_{\cdot 0}$. Moreover it is easy to see that \[X^k = P_{\clg} B^k|_{\clg},\]for all $k \geq 0$. Hence $X^k \raro 0$ strongly as $k \raro 0$. Thus $X \in C_{00}$ and therefore
\begin{equation}\label{e4.3}\mbox{dim~} \cld_X = \mbox{dim~} \cld_{X^*}.\end{equation}
But according to the discussion preceding the statement of Theorem \ref{S-C}, $T|_{\clm}$ admits the decomposition \[T|_{\clm} = T|_{\clm_m} \oplus T|_{\clm_0}\] with respect to the orthogonal splitting  $\clm = \clm_m \oplus \clm_0$, where \[\clm_m = \bigvee_{k \geq 0} (T|_{\clm})^k \clg\]
is the minimal isometric dilation space of $X$ and $\clm_{0} \subseteq \clh_1$. By the uniqueness of the minimal isometric dilation (see \cite{FF}, page 135), $\clm_{m}$ admits a decomposition of the form $$\clm_{m} = \clg \oplus M_+(\varphi \cld_X),$$ where $\varphi$ is a unitary operator from $\cld_X$ to $\mbox{ker} ( (T_{\clm_m})|_{\clm_m \ominus \clg})^*$ and $M_+(\varphi \cld_X) = \varphi \cld_X \oplus T_{\clm_m} (\varphi \cld_X) \oplus \cdots$. Hence $$\clm = (\clg \oplus M_+(\varphi \cld_X)) \oplus \clm_0,$$ which implies that
$$\clh_1 = ( \clh_1 \cap \clm_m) \oplus \clm_0 = M_+(\varphi \cld_X) \oplus \clm_0.$$ Hence, $$\mbox{ker} S^* = \varphi \cld_X \oplus \mbox{ker} (T|_{\clm_0})^*.$$
On the other hand (see \cite{FF}, page 136), there exist a unitary operator $\varphi_*$ from $\cld_{X^*}$ to ker$(T|_{\clm_m})^*$ such that $$\bigvee_{k \geq 0}^{\infty} (T|_{\clm})^k \clg = M_+(\varphi_* \cld_{X^*}),$$ so that
\begin{equation}\label{e4.4} \mbox{ker} (T|_{\clm})^* = \mbox{ker} (T|_{\clm_m})^* \oplus \mbox{ker} (T_{\clm_0})^* = \varphi_* \cld_{X^*} \oplus \mbox{ker} (T|_{\clm_0})^*.\end{equation}
Finally, by (\ref{e4.3}) there exists a unitary operator from $\cld_X$ to $\cld_{X^*}$. Therefore, the relation (\ref{e4.2}) follows from (\ref{e4.4}). This completes the proof of Step 2.

\NI\textsf{Step 3:} Applying Step 2 to the setting in Step 1 we see that \[\mbox{dim~}\mbox{ker~} S'^* = \mbox{dim~}\mbox{ker~} S^*,\]where $S = T_{\clm}$ and $\clm$ is the maximal subspace invariant to $T$ on which the restriction of $T$ is isometric and this completes the proof of the theorem. \qed


\newsection{On the uniqueness of the nilpotent operators}


In this section we will discuss the uniqueness of the nilpotent operators in the representations of c.n.u. contractions with polynomial characteristic functions. We start with the following proposition.

\begin{Proposition}\label{N1}
Let $T$ on $\clh$ be in $\clp_n$ with two different matricial representations
\begin{equation}\label{5.1}
T = \begin{bmatrix} S&*&*\\0&N&*\\0&0&C\end{bmatrix} \; \mbox{and} \; T = \begin{bmatrix} S'&*&*\\0&N'&*\\0&0&C'\end{bmatrix},
\end{equation}
on $\clh = \clh_1 \oplus \clh_0 \oplus \clh_{-1}$ and $\clh = \clh'_1 \oplus \clh'_0 \oplus \clh'_{-1}$, respectively, with $\clm_1 \cap \clh_0 = \{0\} = \clh_0 \cap \clm_{-1}$ and $\clm_1 \cap \clh'_0 = \{0\} = \clh'_0 \cap \clm_{-1}$. Define $$Y = P_{\clh'_0 \oplus \clh'_{-1}}|_{\clh_1 \oplus \clh_0},$$ $$\cll_1 = \clh_1 \ominus (\clh_1 \cap \clh'_1) \;\mbox{and}\; \cll'_{-1} = \clh'_{-1} \ominus (\clh_{-1} \cap \clh'_{-1}).$$ Then:

(i)\; $\mbox{ker~}Y = \clh_1 \cap \clh'_1, \; \mbox{ker~}Y^* = \clh_{-1} \cap \clh'_{-1}$.

(ii)\; $\mbox{ker~} \tilde{Y} = \{0\}$ and $(\mbox{ran}~\tilde{Y})^- = \clh'_0 \oplus \cll'_{-1},$ where $$\tilde{Y} = Y|_{\cll_1 \oplus \clh_0} \colon\ \cll_1 \oplus \clh_0 \rightarrow \clh'_0 \oplus \cll'_{-1}.$$In particular, $\tilde{Y}$ is a quasi-affinity.

(iii)\; $\tilde{Y} \tilde{T} = \tilde{T'} \tilde{Y},$ where $$\tilde{T} = P_{\cll_1 \oplus \clh_0} T|_{\cll_1 \oplus \clh_0}, \; \tilde{T'} = P_{\clh'_0 \oplus \cll'_{-1}} T|_{\clh'_0 \oplus \cll'_{-1}}.$$
\end{Proposition}

\NI \textsf{Proof.} (i) Let $h = h_1 + h_0 \in \clh_1 \oplus \clh_0$ such that $Y h = 0$. Then $h \in \clh'_1$, and hence (since $\clh_1, \clh'_1 \subseteq \clm_1$), $h_0 = h - h_1 \in \clm_1$. Consequently, $h_0 \in \clm_1 \cap \clh_0$ and from the condition $\clm_1 \cap \clh_0 = \{0\}$ we see that $h_0 = 0$. Therefore, $\mbox{ker} Y = \clh_1 \cap \clh'_1$. The second equality follows as above with the observation that $Y^* = P_{\clh_1 \oplus \clh_0}|_{\clh'_0 \oplus \clh'_{-1}}$.

(ii) $\mbox{ker} \tilde{Y} = \{0\}$ follows from the definition of $\tilde{Y}$. The second equality follows from the equalities $$(\mbox{ran}~\tilde{Y})^- = (\mbox{ran}~Y)^- = (\mbox{ker~} Y^*)^{\perp} =  (\clh'_0 \oplus \clh'_{-1}) \ominus (\clh_{-1} \cap \clh'_{-1}) = \clh'_0 \oplus \cll'_{-1}.$$

(iii)
We first notice that, since $$T(\clh_1 \oplus \clh_0) \subseteq \clh_1 \oplus \clh_0, \quad \mbox{and} \quad T^*(\clh'_0 \oplus \clh'_{-1}) \subseteq \clh'_0 \oplus \clh'_{-1},$$ we have $$Y T|_{\clh_1 \oplus \clh_0} = P_{\clh'_0 \oplus \clh'_{-1}} T|_{\clh_1 \oplus \clh_0} = P_{\clh'_0 \oplus \clh'_{-1}} T P_{\clh'_0 \oplus \clh'_{-1}}|_{\clh_1 \oplus \clh_0} = P_{\clh'_0 \oplus \clh'_{-1}} T Y.$$But $$(\mbox{ran}~\tilde{Y})^- = (\mbox{ran}~Y)^- = \clh'_0 \oplus \cll'_{-1},$$therefore $$Y T|_{\clh_1 \oplus \clh_0} = P_{\clh'_0 \oplus \cll'_{-1}} Y T|_{\clh_1 \oplus \clh_0} = P_{\clh'_0 \oplus \cll'_{-1}} T Y = P_{\clh'_0 \oplus \cll'_{-1}} T|_{\clh'_0 \oplus \cll'_{-1}} Y,$$that is $$Y T|_{\clh_1 \oplus \clh_0}  = \tilde{T}' Y.$$Moreover, we have $$Y T|_{\clh_1 \oplus \clh_0} = Y (I - P_{\mbox{ker}\, Y}) T|_{\clh_1 \oplus \clh_0} =  Y (I - P_{\mbox{ker}\, Y}) T  (I - P_{\mbox{ker}\, Y})|_{\clh_1 \oplus \clh_0} = \tilde{Y} \tilde{T} P_{\cll_1 \oplus \clh_0}|_{\clh_1 \oplus \clh_0},$$ and thus $$\tilde{Y} \tilde{T} = Y T|_{\cll_1 \oplus \clh_0} = \tilde{T}' Y|_{\cll_1 \oplus \clh_0} = \tilde{T}' \tilde{Y}.$$ This finishes the proof of the proposition. \qed

\vspace{0.1in}

Let
\begin{equation}\label{a}
\begin{bmatrix} S&*&*\\0&N&*\\0&0&C\end{bmatrix}\quad \mbox{and} \quad \begin{bmatrix} {S_*}&*&*\\0&{N_*}&*\\0&0&{C_*}\end{bmatrix}
\end{equation}
be the matrix representations of $T$ on $\clh = \clh_1 \oplus \clh_0 \oplus \clh_{-1}$ and on $\clh = \clh_{1*} \oplus \clh_{0*} \oplus \clh_{-1*}$ considered in Theorem \ref{prop1} and Remark \ref{rem3.5}, respectively. As observed in Section 3 we have that
$$\clh_{-1} = \clm_{-1} \; \; \mbox{and} \; \; \clh_{1*} = \clm_1,$$thus $$ \clh_0 \cap \clm_{-1} = \{0\} \; \; \mbox{and}\;\; \clm_1 \cap \clh_{0*} = \{0\}.$$However, we do not know if $$\clm_1 \cap \clh_0 = \{0\} \;\; \mbox{and} \;\; \clh_{0*} \cap \clm_{-1} = \{0\}.$$Therefore, the representations in (\ref{a}) may not satisfy the hypotheses of Proposition \ref{N1}. However, note that by defining \begin{equation}\label{b} \clh_{10} = \clh_1 \oplus ( \clm_1 \cap \clh_0), \; \clh_{00} = \clh_0 \ominus (\clm_1 \cap \clh_0), \; \clh_{-10} = \clh_{-1}\end{equation} and \begin{equation}\label{c}\clh_{1*0} = \clh_{1*}, \; \clh_{0*0} = \clh_{0*} \ominus (\clh_{0*} \cap \clm_{-1}),\; \clh_{-1*0} = \clh_{-1*} \oplus (\clh_{0*} \cap \clm_{-1}),\end{equation}the representations
\begin{equation}\label{d}
\begin{bmatrix} S_0&*&*\\0&N_0&*\\0&0&C_0\end{bmatrix}\quad \mbox{and} \quad \begin{bmatrix} {S_{*0}}&*&*\\0&{N_{*0}}&*\\0&0&{C_{*0}}\end{bmatrix},
\end{equation}
on $\clh = \clh_{10} \oplus \clh_{00} \oplus \clh_{-10}$ and $\clh = \clh_{1*0} \oplus \clh_{0*0} \oplus \clh_{-1*0}$, respectively satisfy those hypotheses.

\begin{Theorem}\label{T5.2}
The "central" nilpotent operators $N_0$ and $N_{*0}$ in the representations (\ref{d}) are quasi-similar.
\end{Theorem}
\NI \textsf{Proof.} Let the representations in (\ref{d}) play the role of the representations in (\ref{N1}), that is, with $\clh_1 = \clh_{10}, \clh_0 = \clh_{00}, \clh_{-1} = \clh_{-10}$ and $\clh'_1 = \clh_{1*}, \clh'_0 = \clh_{0*0}, \clh'_{-1} = \clh_{-1*0}$, respectively. Then due to (\ref{b}) we have \begin{equation*}\cll_{1} = \clh_{10} \ominus (\clh_{10} \cap \clh_{1*}) = \{0\} \; \mbox{and} \; \cll'_{-1} = \clh_{-1*0} \ominus ( \clh_{-1*0} \cap \clh_{-1}) = \{0\}.\end{equation*}
Therefore, $\tilde{T} = N_0, \, \tilde{T}' = N_{*0}$ and Proposition \ref{N1} (ii) and (iii) shows that $N_0$ is a quasi-affine transform of $N_{*0}$, or, in the standard notation, $N_0 \prec N_{*0}$. Since these operators are of class $C_0$, it follows that $N_0$ and $N_{*0}$ are quasi-similar (see Proposition 5.1 in \cite{B}). This concludes the proof.\qed

To state the next result we recall that  $T_1 \in \cll(\clh_1)$ is said to be injected in $T_2 \in \cll(\clh_2)$ which is denoted by $T_1 \prec^i T_2$ if there is an injective $X \in \cll(\clh_1, \clh_2)$ such that $X T_1 = T_2 X$.

\begin{Theorem}\label{5.3}
Let \begin{equation}\label{e}
T = \begin{bmatrix} S&*&*\\0&N&*\\0&0&C\end{bmatrix},
\end{equation}
be any matrix representation on $\clh = \clh_1 \oplus \clh_0 \oplus \clh_{-1}$ such that $\clm_1 \cap \clh_0 = \{0\} = \clh_0 \cap \clm_{-1}$. Then $$N_{*0} \prec^i N\prec^i N_*.$$
\end{Theorem}
\NI \textsf{Proof.} We will apply Proposition \ref{N1} with the representation (\ref{e}) in the role of the first representation in (\ref{5.1}) and with the second representation in (\ref{5.1}) replaced by the second representation in (\ref{d}). Then we have $$\cll_1 = \clh_1 \ominus ( \clh_1 \cap \clh_{1*0}) = \clh_1 \ominus ( \clh_1 \cap \clm_1) = \{0\}.$$Consequently, by applying Proposition \ref{N1}, (ii) and (iii) we obtain $$N = \tilde{T} \prec \tilde{T}' = P_{\clh_{0*0} \oplus \cll'_{-1}} T|_{\clh_{0*0} \oplus \cll'_{-1}},$$ where $$\cll'_{-1} = \clh_{-1*0} \ominus (\clh_{-1} \cap \clh_{-1*0}),$$ and we infer that $\tilde{T}'$ is a nilpotent operator and then, since both $\tilde{T}$ and $\tilde{T}'$ are of class $C_0$, that $\tilde{T} \sim \tilde{T}'$. Observe now that $$\clh_{-1*0} = (\clh_{0*} \oplus \clh_{-1*}) \cap \clm_{-1},$$ and $$\clh_{-1*0} \cap \clh_{-1} = (\clh_{0*} \oplus \clh_{-1*}) \cap \clh_{-1} = (\clh_{0*} \cap \clh_1) \oplus \clh_{-1*},$$by the minimality property of $\clh_{-1*}$ (see Remark \ref{rem3.6}). Therefore,
\begin{equation*}
\begin{split}\cll'_{-1}  & = [(\clh_{0*} \cap \clm_{-1}) \oplus \clh_{-1*}] \ominus [(\clh_{0*} \cap \clh_{-1}) \oplus \clh_{-1*}]\\ & = (\clh_{0*} \cap \clm_{-1}) \ominus ( \clh_{0*} \cap \clh_{-1}),
\end{split}
\end{equation*}
and
\begin{equation*}
\begin{split}
\clh_{0*0} \oplus \cll'_{-1} & = [\clh_{0*} \ominus ( \clh_{0*} \cap \clm_{-1})] \oplus [ (\clh_{0*} \cap \clm_{-1}) \ominus (\clh_{0*} \cap \clh_{-1})] \\ & = \clh_{0*} \ominus ( \clh_{0*} \cap \clh_{-1}),
\end{split}
\end{equation*}
and thus we see that
\begin{equation*}
\begin{split}
\tilde{T}' & = P_{\clh_{0*} \ominus ( \clh_{0*} \cap \clh_{-1})} T|_{\clh_{0*} \ominus ( \clh_{0*} \cap \clh_{-1})} \\ & = P_{\clh_{0*} \ominus ( \clh_{0*} \cap \clh_{-1})} N_*|_{\clh_{0*} \ominus ( \clh_{0*} \cap \clh_{-1})}.
\end{split}
\end{equation*}
But $N_*^* (\clh_{0*} \cap \clh_{-1}) = P_{\clh_{0*}} T^*(\clh_{0*} \cap \clh_{-1})$ and for  $h \in \clh_{0*} \cap \clh_{-1}$ we have $$T^* h \in (\clh_{0*} \oplus \clh_{-1*}) \cap \clh_{-1} = (\clh_{0*} \cap \clh_{-1}) \oplus \clh_{-1*},$$where $$N_* h = P_{\clh_{0*}} T^* h \in \clh_{0*} \cap \clh_1;$$ consequently we have $$N_* ( \clh_{0*} \cap \clh_{-1}) \subseteq \clh_{0*} \cap \clh_{-1}.$$ Therefore, we infer \begin{equation*} N_* (\clh_{0*} \ominus ( \clh_{0*} \cap \clh_{-1})) \subseteq \clh_{0*} \ominus ( \clh_{0*} \cap \clh_{-1}).\end{equation*} Now it follows that \begin{equation}\label{f}N \sim \tilde{T}' = N_*|_{\clh_{0*} \ominus ( \clh_{0*} \cap \clh_{-1})}.\end{equation} In particular, (\ref{f}) shows that $N \prec^i  N_*$. Also since $\clh_{0*0} = \clh_{0*} \ominus (\clh_{0*} \cap \clm_{-1}),$ we have that $\clh_{0*0} \subseteq \clh_{0*} \ominus (\clh_{0*} \cap \clh_{-1})$ and $$N_*|_{\clh_{0*0}} = \tilde{T'}|_{\clh_{0*0}}.$$ Consequently, $$N_{*0} = \tilde{T}'|_{\clh_{*0} \ominus ( \clh_{*0} \cap \clh_{-1})},$$ and the relation (\ref{f}) implies that $N_{*0} \prec^i N$. \qed

The following remark follows by applying Theorem \ref{5.3} to $T^*$.
\begin{Remark}\label{5.4}
Under the conditions of Theorem \ref{5.3} one obtains $$N^*_0 \prec^i N^* \prec^i N_1^*,$$where $N_1$ denotes the canonical nilpotent operator as in Theorem \ref{prop1}. Note that since $N_0 \sim N_{*0}$ we have $N_0 \prec^i N$ and $N^*_{*0} \prec^i N^*_1.$
\end{Remark}

Let $T$ be a c.n.u. contraction with polynomial characteristic function of degree $n$ and let
\begin{equation*}
T = \begin{bmatrix} S&*&*\\0&N_{um}&*\\0&0&C\end{bmatrix}.
\end{equation*}
Then $N_{um}$ is said to be \textit{ultra-minimal} if $N_{um} \prec^i N$ for any other nilpotent operator $N$ of the upper triangular matrix representations of $T$.

By Theorem \ref{5.3}, the canonical $N_{*0}$ is ultra-minimal. If $N$ is any other ultra-minimal nilpotent operator in a upper triangular matrix representation of $T$, then we will have $$N' \prec^i N_{*0} \prec^i N',$$
so there exists injective bounded linear operators $X \colon\ \clh'_0 \raro \clh_{0*0}$ and $Y \colon\ \clh_{0*0} \raro \clh'_0$ such that
\begin{equation}\label{5.8}
X N' = N_{*0} X,
\end{equation}
and
\begin{equation}
Y N_{*0} = N' Y.
\end{equation}
Thus $$N_{*0} X Y = X Y N_{*0},$$ hence $$N_{*0} \prec N_{*0}|_{(\mbox{ran} X Y)^-} = N_{*0}|_{(\mbox{ran} X)^-}.$$ But since $N_{*0}$ and $N_{*0}|_{(\mbox{ran} X)^-}$ are $C_0$-operators we have (see \cite{B}) $$N_{*0} \sim N_{*0}|_{(\mbox{ran} X)^-}.$$Since (\ref{5.8}) yields also $$N' \prec N_{*0}|_{(\mbox{ran} X)^-},$$ we have $$N' \sim N_{*0}|_{(\mbox{ran} X)^-} \sim N_{*0}.$$ Consequently all ultra-minimal nilpotent operators are quasi-similar and hence we have the following theorem.

\begin{Theorem}\label{5.5}
The family of all ultra-minimal nilpotent operators in the representations of type (\ref{1}) is a quasi-similar equivalence class (containing $N_0$ and $N_{*0}$) which moreover, is a unitary invariant of $T$.
\end{Theorem}

\begin{Remark}\label{5.6}
Recall that in the diagonal of the representations (\ref{1}) of $T$, dim $(\mbox{ker}\,S^*)$ and dim $(\mbox{ker}\,C)$ are unitary invariants of $T$. Obviously, they can be replaced by the unitary equivalent classes of $S$ and $C$. Thus we still have the following open problem. Can the quasi-similar class considered in Theorem \ref{5.5} be replaced with the unitary equivalence class?
\end{Remark}

\newsection{Degenerate forms}

Let $T \in \cll(\clh)$ be a c.n.u. contraction with a operator-valued polynomial characteristic function $\Theta_T$. Then, as was proved in \cite{T}, \cite{WP}, \cite{CM} and \cite{BM} that the representation (\ref{1}) take the degenerate form

\begin{equation*}
\begin{bmatrix} S&*\\0&C\end{bmatrix},
\end{equation*}
if and only if $$\Theta_T(z) \equiv \Theta_T(0), \quad (z \in \mathbb{D}).$$ For the other degenerate forms of (\ref{1}) we have the following.

\begin{Proposition}
Let $T \in \cll(\clh)$ be a c.n.u. contraction with an operator-valued polynomial characteristic function. Then (\ref{1}) has the degenerate form

(i) $T = \begin{bmatrix} S& *\\ 0& N\end{bmatrix}$ if and only if $T \in C_{\cdot \,0}$.

(ii) $T = \begin{bmatrix} N& *\\ 0& C\end{bmatrix}$ if and only if $T \in C_{0\, \cdot}$.

(iii) $T = \begin{bmatrix} S\end{bmatrix}$ if and only if $T \in C_{1\,0}$, that is if and only if $\Theta_T\equiv 0 \in \cll(\{0\}, \cle_*)$ for some Hilbert space $\cle_* \neq \{0\}$.

(iv) $T = \begin{bmatrix} C\end{bmatrix}$ if and only if $T \in C_{0\,1}$, that is if and only if $\Theta_T\equiv 0 \in \cll(\cle, \{0\})$ for some Hilbert space $\cle \neq \{0\}$.

(v) $T = \begin{bmatrix} N\end{bmatrix}$ if and only if $T \in C_{0\,0}$.
\end{Proposition}

\NI \textsf{Proof.} (i) If $T = \begin{bmatrix} S& *\\ 0& N\end{bmatrix}$, then $$T^{*n} = \begin{bmatrix} S^{*n}&0\\X& 0\end{bmatrix},$$for some $X \in \cll(\clh_0, \clh_1)$ and hence $$(T^{*n})^k (h_1 \oplus h_0) = S^{* n k} h_1 \oplus  X S^{*n(k-1)} h_0, \quad \quad (h_1 \oplus h_0 \in \clh_1 \oplus \clh_0).$$ Therefore,
$$
\mbox{lim}_{k \raro \infty} \|T^{*nk}(h_1 \oplus h_0)\| = \mbox{lim}_{k \raro \infty} \|(T^{*n})^k (h_1 \oplus h_0)\|^2 = \mbox{lim}_{k \raro \infty} (\|S^{*nk}h_1\|^2 + \|X S^{*n(k-1)} h_0\|^2) = 0.
$$
So $T \in C_{\cdot \, 0}$. Conversely if $$T = \begin{bmatrix} S&*&*\\0&N&*\\0&0&C\end{bmatrix} \in C_{\cdot\,0},$$ then for any $h_{-1} \in \clh_{-1}$, we have that $$\|h_{-1}\| = \|C^{*k} h_{-1}\| = \|T^{*k}(0 \oplus 0 \oplus h_{-1})\| \raro 0, \quad \mbox{as} k \raro \infty.$$Thus we obtain that $\clh_{-1} = \{0\}$ and $T$ takes the degenerate form $T = \begin{bmatrix}S&*\\0&N\end{bmatrix}$. Finally, $T \in C_{\cdot\,0}$ if and only if $\Theta_T$ is inner (see \cite{NF}, Chapter VI).

\NI (ii) follows from (i) applied to $T^*$.

\NI (iii) follows from the fact that the characteristic function of any shift coincides with $0 \in \cll(\cld_S, \cld_{S^*}) = \cll(\{0\}, \mbox{ker} S^*)$.

\NI (iv) Apply (iii) to $T^*$.

\NI (v) $T = \begin{bmatrix} N\end{bmatrix}$ if and only if the representation is a degenerate representations of both forms in (i) and (ii), hence if and only if $T \in C_{0\,0}$.

\newsection{The case of monomial characteristic functions}

In this section we study the the c.n.u. contraction $T$ such that $\Theta_T(z)$ is a non-constant monomial (for the case of constant monomials, see Section 6). This particular case, although quite concrete, allow us to show that the nilpotent central entries in the two canonical representations may not be quasi-similar. For this purpose, we explicitly (in terms of $\Theta_T$) calculate the two canonical triangular representations of $T$. First, we need the following simple fact, which, for instance, can be obtained by using (\ref{ch-def}).

\vspace{0.1in}

\begin{Lemma}\label{unitary-part}
Let $m \geq 1$ and $U$ be an unitary operator in $\cll(\cle, \cle_*)$ and $\Theta(z) = U z^m$. Then the c.n.u contraction $J_m(\cle_*)$ for which $\Theta_{J_m(\cle_*)}(z)$ coincides with $\Theta(z)$ is a nilpotent operator of order $m$. Moreover,
\begin{equation*}
J_m(\cle_*) = \begin{bmatrix} 0&I_{\cle_*}&0&\cdots &0&0&0\\0&0&I_{\cle_*}& \cdots & 0&0&0\\ 0& 0&0& \cdots & 0&0&0\\ \vdots & \vdots &\vdots& \cdots& \vdots &\vdots &\vdots \\ 0&0&0& \cdots & 0&I_{\cle_*}&0\\0 & 0 & 0&\cdots& 0&0&I_{\cle_*}\\0&0&0&\cdots &0&0&0\end{bmatrix}.
\end{equation*}
\end{Lemma}

\vspace{0.2in}

Let $A$ in $\cll(\clm, \cln)$ be a non-zero contraction. Then it is easy to see that $A|_{{\cld_A}^{\perp}} \colon\ {{\cld_A}^{\perp}} \raro {{\cld_{A^*}}^{\perp}}$ is an unitary operator and $A|_{\cld_A} \colon\ \cld_A \raro \cld_{A^*}$ is pure contraction, that is, $\|A h\| < \|h\|$ for all $h \neq 0$ in $\cld_A$. Now, if $\Theta(z) = A z^m$ is a given non-constant monomial where $A$ in $\cll(\clm, \cln)$ is a contraction then $\Theta(z) = A|_{\cld_A} z^m \oplus A|_{{\cld_A}^{\perp}} z^m$. In what follows, we will assume that $m \geq 1$. In view of Lemma \ref{unitary-part} we have that, $T$ is unitarily equivalent to $T_p \oplus J_m(\cle_*)$ where $T_p$ is the c.n.u. contraction whose characteristic function coincides with $A|_{\cld_A} z^m$. Now for a given pure contraction $A$ in $\cll(\clm, \cln)$, we will construct the two canonical upper triangular representations  $T_A$ and $\tilde{T}_A$ with the characteristic function (coinciding with) $\Theta (z) = A z^m$ for $z$ in $\mathbb{D}$.

\NI Define a Hilbert space $\clh = \oplus_{- \infty}^{\infty} \cll_n$, where
\begin{equation*}
 \cll_n = \left \{
\begin{array}{cl}
\cln & \mbox{if } n \geq 0\\
\clm & \mbox{if } n \leq - 1\\
\end{array}
\right.
\end{equation*}
For convenience, the element $\bm{h}$ of $\clh$ will be denoted by $\bm{h} = \cdots \oplus h_1 \oplus h_0 \oplus h_{-1} \oplus \cdots$ where the subscript refer that $h_j \in \cll_j$ for all $j \in \mathbb{Z}$. Also, we will denote the subspace $\cdots \oplus \{0\} \oplus \{0\} \oplus \cll_j \oplus \{0\} \oplus \{0\} \oplus \cdots$ of $\clh$ simply by $\cll_j$.

We define a bounded linear operator $T_A$ acting on $\clh$ by

\begin{equation*}
T_A|_{\cll_n} = \left \{
\begin{array}{cl}
D_{A^*} \colon\ \cll_{m-1} \raro \cll_{m} & \mbox{if } n = m-1\\
- A\colon\ \cll_{-1} \raro \cll_{m} & \mbox{if } n = -1\\
I \colon\ \cll_n \raro \cll_{n+1} & \mbox{if } n \neq -1, m-1.
\end{array}
\right.
\end{equation*}

Then $T_A$ has the following matrix representation

\begin{equation}\label{bigT_A}
T_A = \left[
  \begin{array}{ c c c c c|cc c c c |c c c c c}
\vdots& \cdots&\vdots &\vdots&\vdots  &\vdots& \vdots &\cdots&\vdots&\vdots &\vdots& \vdots&\vdots&\cdots&\vdots\\
\vdots&\ddots&I_{\cln}&0&0 &0& 0&\cdots&0&0 &0& 0&0&\cdots&\vdots\\
\vdots&\cdots&\bm{0}&I_{\cln}&0 &0& 0&\cdots&0&0 &0& 0&0&\cdots&\vdots\\
\vdots&\cdots&0&\bm{0}&I_{\cln} &0& 0&\cdots&0&0 &0& 0&0&\cdots&\vdots\\
\vdots&\cdots&0&0&\bm{0} &D_{A^*}& 0&\cdots&0&0 &- A& 0&0&\cdots&\vdots\\ \hline

\vdots&\cdots&0&0&0 &\bm{0}& I_{\cln}&\cdots&0&0 &0& 0&0&\cdots&\vdots\\
\vdots&\cdots&0&0&0 &0& \bm{0}&\cdots&0&0 &0& 0&0&\cdots&\vdots\\
\vdots&\cdots&\vdots&\vdots&\vdots &\vdots& \vdots&\cdots&\vdots&\vdots &\vdots& \vdots&\vdots&\cdots&\vdots\\
\vdots&\cdots&0&0&0 &0& 0&\cdots&\bm{0}&I_{\cln}&0& 0&0&\cdots&\vdots\\

\vdots&\cdots&0&0&0 &0& 0&\cdots&0&\bm{0} &0& 0&0&\cdots&\vdots\\ \hline

\vdots&\cdots&0&0&0 &0& 0&\cdots&0&0 &\bm{0}& I_{\clm}&0&\cdots&\vdots\\
\vdots&\cdots&0&0&0 &0& 0&\cdots&0&0 &0& \bm{0}&I_{\clm}&\cdots&\vdots\\
\vdots&\cdots&0&0&0 &0& 0&\cdots&0 &0& 0&0&\bm{0}&\cdots&\vdots\\
\vdots&\cdots&0&0&0 &0& 0&\cdots&0 &0& 0&0&0&\ddots&\vdots\\
\vdots&\cdots&\vdots&\vdots&\vdots &\vdots& \vdots&\cdots&\vdots&\vdots &\vdots& \vdots&\vdots&\cdots&\vdots\\

\end{array} \right],
\end{equation}
where the central block matrix is of order $m$.

To compute the defect operators of $T_A$, we first observe that

\begin{equation*}
T^*_A|_{\cll_n} = \left \{
\begin{array}{cl}
0 \colon\ \cll_0 \raro \cll_{-1} & \mbox{if } n = 0\\
D_{A^*} \oplus -A^* \colon\ \cll_{m} \raro  \cll_{m-1} \oplus \cll_{-1} & \mbox{if } n = m\\
I \colon\ \cll_n \raro \cll_{n-1} & \mbox{if } n \neq 0, m,
\end{array}
\right.
\end{equation*}

so

\begin{equation*}
T_A T^*_A|_{\cll_n} = \left \{
\begin{array}{cl}
0 \colon\ \cll_0 \raro \cll_0 & \mbox{if } n = 0\\
I \colon\ \cll_n \raro \cll_n & \mbox{if}\; n \neq 0,
\end{array}
\right.
\end{equation*}

and

\begin{equation*}
T^*_A T_A|_{\cll_n} = \left \{
\begin{array}{cl}
D_{A^*}^2 \oplus -A^* D_{A^*} \colon\ \cll_{m-1} \raro \cll_{m-1} \oplus \cll_{-1} & \mbox{if } n = m-1\\
- D_{A^*} A \oplus A^* A \colon\ \cll_{-1} \raro \cll_{m-1} \oplus \cll_{-1} & \mbox{if } n = -1\\
I \colon\ \cll_n \raro \cll_n & \mbox{if } n \neq m-1, -1
\end{array}.
\right.
\end{equation*}

Thus, the defect operators are given by

\begin{equation}\label{71a}
D^2_{T^*_A}|_{\cll_n} = \left \{
\begin{array}{cl}
I \colon\ \cll_0 \raro \cll_0 & \mbox{if } n = 0\\
0 \colon\ \cll_n \raro \cll_n  & \mbox{if } n \neq 0,
\end{array}
\right.
\end{equation}

\begin{equation}\label{71b}
D^2_{T_A}|_{\cll_{m-1} \oplus \cll_{-1}} =
\begin{bmatrix}
A A^* & A D_{A}\\D_{A} A^* & D^2_{A}\end{bmatrix} \colon\ \cll_{m-1} \oplus \cll_{-1} \raro \cll_{m-1} \oplus \cll_{-1},
\end{equation}

and $$D^2_{T_A}|_{\cll_n} = 0 \colon\ \cll_n \raro \cll_n, \;\;\; \mbox{if} \; n \neq m-1, -1.$$

From (\ref{71a}) we infer that $T_A$ is a partial isometry and therefore, the defect spaces of $T_A$ can be expressed as
\begin{equation*}
 \begin{split}\cld_{T_A}  = & \cdots \{0\} \oplus \{ A (A^* h_{m-1} + D_A h_{-1})  \oplus 0 \oplus \cdots \oplus  0 \oplus D_A (A^* h_{m-1} + D_A h_{-1}) \colon\ \\ & \;\;\;\; h_{m-1} \in \clh_{m-1}, h_{-1} \in \clh_{-1}\}^- \oplus \{0\} \oplus \cdots\\ = & \{A v \oplus D_A v \in \cll_{m-1} \oplus \cll_{-1} \colon\ v \in \clm\},
\end{split}
\end{equation*}
and
$$\cld_{T^*_A} =  \cdots \oplus \{0\} \oplus \cll_{0} \oplus \{0\} \cdots.$$

\NI Also since $A$ is a pure contraction, we have the following.

\begin{Lemma}
The operator $T_A$ constructed as above is c.n.u.
\end{Lemma}
\NI\textsf{Proof.}  Let $\clr$ be a $T_A$-reducing closed subspace of $\clh$ such that $T_A|_{\clr}$ is unitary. We need to prove that $\clr = \{\bm{0}\}$.

\NI Consider a vector $\bm{h} = \cdots \oplus h_{m+1} \oplus h_{m} \oplus h_{m-1} \oplus h_{m-2} \oplus \cdots \oplus h_{1} \oplus h_{0} \oplus h_{-1} \oplus h_{-2} \oplus \cdots$ from $\clr$. Since $\clr$ is $T_A$-reducing and $T_A|_{\clr}$ is unitary, $D^2_{T^*_A} \bm{h} = \bm{0}$, which yields $h_{0} = 0$. On the other hand, it is easy to see that
\begin{equation*}
 \begin{split}T^m_A \bm{h} = & \cdots \oplus (h_m)_{2m+1} \oplus (D_{A^*} h_{m-1} - A h_{-1})_{2m} \oplus \cdots \oplus (D_{A^*} h_0 - A h_{-m})_m \oplus \\& (0)_{m-1} \oplus \cdots \oplus (0)_0 \oplus (h_{-m-1})_{-1} \oplus (h_{-m-2})_{-2} \oplus \cdots.
\end{split}
\end{equation*}
Since $T_A^m \clr = \clr$, it follows that for any $\bm{h}$ in $\clr$, $$h_{m-1} = h_{m-2} = \cdots= h_1 = h_0 = 0.$$
Hence
\begin{equation}\label{7.3}
 \begin{split}T^{m+1}_A \bm{h} = & \cdots \oplus (h_m)_{2m+2} \oplus ( - A h_{-1})_{2m+1} \oplus \cdots \oplus ( - A h_{-m})_{m+1} \oplus  (- A h_{-m-1})_m \\& \oplus  (0)_{m-1} \oplus \cdots \oplus (0)_0 \oplus (h_{-m-2})_{-1} \oplus (h_{-m-3})_{-2} \oplus \cdots.
\end{split}
\end{equation}
 Consequently, by the facts that $T_A|_{\clr}$ is unitary and $A$ is a strict contraction, (\ref{7.3}) yields $h_{-1} = h_{-2} = \cdots = h_{-m-1} = 0$. By repeated application of $T_A$ on $T^{m+1}_A \bm{h}$ we conclude that $$h_{-k} = 0, \quad \mbox{for all}\; k\geq 1.$$ Therefore, we obtain for all $\bm{h}$ in $\clr$, $$\bm{h} = \cdots \oplus h_{m+1} \oplus h_m \oplus 0 \oplus 0 \oplus \cdots \oplus 0 \oplus 0 \oplus 0 \oplus \cdots.$$
But now the $(m-1)$-th component, $(D_{A^*} h_m)_{m-1}$, of
\begin{equation*}
\begin{split}
T^*_A \bm{h} = &\cdots \oplus (h_{m+1})_m \oplus (D_{A^*} h_m)_{m-1} \oplus \cdots \oplus (h_2)_1 \oplus (h_1)_0 \oplus (-A^* h_m)_{-1} \oplus (h_{-1})_{-2} \\& \oplus \cdots \oplus (h_{-m})_{-m-1} \oplus (0)_{-m-2} \oplus \cdots,
\end{split}
\end{equation*}
must be zero. Since $\mbox{ker}\, D_{A^*} = \{0\}$ then $h_m$ must be the zero vector. Moreover, using the fact $T_A^* \clr = \clr$, we infer that
$$\bm{h} = \cdots \oplus h_{m+2} \oplus h_{m+1} \oplus 0 \oplus 0 \oplus 0 \cdots.$$
Finally, repeating the above argument replacing the role of $T_A$ by $T_A^*$, we obtain $h_k = 0$ for all $k > m$. Thus $\bm{h} = \bm{0}$. \qed

\vspace{0.2in}

\begin{Theorem}\label{Cz}
Let $A \colon\ \clm \raro \cln$ be a pure contraction (that is, $\|A h\|_{\cln} < \|h\|_{\clm}$ for all $h \in \clm$ and $h \neq 0$). Then $T_A$ is a c.n.u. partial isometry and the characteristic function of $T_A$ coincides with $\Theta(z) = A z^m$. Moreover, the upper triangular representation of $T_A$ in (\ref{bigT_A}) coincides with the canonical representation introduced in Theorem \ref{prop1}.
\end{Theorem}

\NI\textsf{Proof.} Let $u =  A v \oplus D_A v \in \cll_{m-1} \oplus \cll_{-1}$ where $v \in \clm$. Now $$ - T_A u = - D_{A^*} A v + A D_A v = 0 \in \cll_m.$$ First, we consider the case when $m > 1$. Since $T_A$ is a partial isometry, we have $$D_{T_A} u = D^2_{T_A} u = u,$$ so that $$D_{T^*_A} D_{T_A} u = D_{T^*_A} u = P_{\cll_0} u = 0.$$ Notice that for $0 \leq j < m-1$,
$$T^{*j}_A u = \cdots \oplus (0)_{m-j} \oplus (Av)_{m-j-1} \oplus (0)_{m-j-2} \oplus \cdots \oplus (0)_{-j} \oplus (D_A v)_{-j-1} \oplus (0)_{-j-2} \oplus \cdots.$$
Consequently for $j = m-1$,
$$T^{*( m-1)}_A u = \cdots \oplus (0)_1 \oplus (Av)_{0} \oplus (0)_{-1} \oplus \cdots \oplus (0)_{-m+1} \oplus (D_A v)_{-m} \oplus (0)_{-m-1} \oplus \cdots,$$
so that $$D_{T^*} T^{*(m-1)} D_T u = P_{\cll_0} T^{*(m-1)} u = A v \in \cll_0,$$
and for $j > m-1$,
\begin{equation*}\begin{split}
D_{T^*_A} T^{*j} D_{T_A} u  & = P_{\cll_0} T^{*j} u \\& =  P_{\cll_0} (\cdots \oplus (0)_{m-j} \oplus (Av)_{m-j-1} \oplus (0)_{m-j-2} \oplus \cdots \oplus \\& \quad \quad (0)_{-j} \oplus (D_A v)_{-j-1} \oplus (0)_{-j-2} \oplus \cdots) \\ &= 0.
\end{split}
\end{equation*}
To finish the proof, we will show that $\tau_* \Theta_{T_A}(z) = \Theta(z) \tau^* $ for all $z \in \mathbb{D}$ for some unitary operators $\tau \colon\ \clm \raro \cld_{T_A}$ and $\tau_* \colon\ \cld_{T^*_A} \raro \cln$. For this purpose, we define $$\tau_* = \cdots \oplus 0 \oplus I_{\cll_{0}} \oplus 0 \oplus \cdots \colon\ \cld_{T^*_A} \raro \cln$$
and
\begin{equation*} \tau = \begin{bmatrix} A\\D_A\end{bmatrix}  \colon\ \clm \raro \cld_{T_A},
\end{equation*}
defined by $$\tau v = Av \oplus D_A v = u, \quad v \in \clm.$$Hence
\begin{equation}\label{7.6a}\Theta_{T_A}(z) u = z^m A v = \tau_*^* \Theta(z) \tau^* u,
\end{equation}
for all $z$ in $\mathbb{D}$.

\NI In the remaining case $m = 1$, taking $u =  A v \oplus D_A v \in \cll_{0} \oplus \cll_{-1}$ we have $$D_{T^*_A} D_{T_A} u = D_{T^*_A} u = P_{\cll_0} u = A v.$$Moreover, for all $k \geq 1$, $$D_{T^*} T^{*k} D_T u = 0,$$ and the rest of the proof follows in the similar way as in the last part of the proof of (\ref{7.6a}) for the $m >1$ case.

\NI Finally, recalling that $A$ is a pure contraction, the last part of the theorem follows by inspection.
\qed

\vspace{0.1in}

In the following, we relate the class of operators with non-constant monomial characteristic functions and the canonical upper triangular representations considered in Remark \ref{rem3.5}. First, given a pure contraction $A \in \cll(\clm, \cln)$, define the bounded linear operator $\tilde{T}_A$ on $\clh = \oplus_{-\infty}^{\infty} \cll_n$  where
\begin{equation*}
 \cll_n = \left \{
\begin{array}{cl}
\cln & \mbox{if } n \geq 0\\
\clm & \mbox{if } n \leq - 1\\
\end{array}
\right.
\end{equation*}
and

\begin{equation*}
\tilde{T}_A|_{\cll_n} = \left \{
\begin{array}{cl}
-A \oplus D_A \colon\ \cll_{-m} \raro \cll_1 \oplus \cll_{-m+1} & \mbox{if } n = - m\\
0\colon\ \cll_{0} \raro \cll_{1} & \mbox{if } n = -1\\
I \colon\ \cll_n \raro \cll_{n+1} & \mbox{if } n \neq 0, - m.
\end{array}
\right.
\end{equation*}

Then $\tilde{T}_A$ admits the following matrix representation

\begin{equation}\label{tildeT_A}
\tilde{T}_A = \left[
   \begin{array}{c c c c c|cc c c c |c c c c c}
\vdots&\cdots&\vdots &\vdots&\vdots  &\vdots& \vdots &\cdots&\vdots&\vdots &\vdots& \vdots&\vdots&\cdots&\vdots\\
\vdots&\ddots&I_{\cln}&0&0 &0& 0&\cdots&0&0 &0& 0&0&\cdots&\vdots\\
\vdots&\cdots&\bm{0}&I_{\cln}&0 &0& 0&\cdots&0&0 &0& 0&0&\cdots&\vdots\\
\vdots&\cdots&0&\bm{0}&I_{\cln} &0& 0&\cdots&0&0 &0& 0&0&\cdots&\vdots\\
\vdots&\cdots&0&0&\bm{0} &0& 0&\cdots&0&0 &- A& 0&0&\cdots&\vdots\\ \hline

\vdots&\cdots&0&0&0 &\bm{0}& I_{\clm}&\cdots&0&0 &0& 0&0&\cdots&\vdots\\
\vdots&\cdots&0&0&0 &0& \bm{0}&\cdots&0&0 &0& 0&0&\cdots&\vdots\\
\vdots&\cdots&\vdots&\vdots&\vdots &\vdots& \vdots&\cdots&\vdots&\vdots &\vdots& \vdots&\vdots&\cdots&\vdots\\
\vdots&\cdots&0&0&0 &0& 0&\cdots&\bm{0}&I_{\clm}&0& 0&0&\cdots&\vdots\\

\vdots&\cdots&0&0&0 &0& 0&\cdots&0&\bm{0} &D_A& 0&0&\cdots&\vdots\\ \hline

\vdots&\cdots&0&0&0 &0& 0&\cdots&0&0 &\bm{0}& I_{\clm}&0&\cdots&\vdots\\
\vdots&\cdots&0&0&0 &0& 0&\cdots&0&0 &0& \bm{0}&I_{\clm}&\cdots&\vdots\\
\vdots&\cdots&0&0&0 &0& 0&\cdots&0 &0& 0&0&\bm{0}&\cdots&\vdots\\
\vdots&\cdots&0&0&0 &0& 0&\cdots&0 &0& 0&0&0&\ddots&\vdots\\
\vdots&\cdots&\vdots&\vdots&\vdots &\vdots& \vdots&\cdots&\vdots&\vdots &\vdots& \vdots&\vdots&\cdots&\vdots\\

\end{array} \right].
\end{equation}

The proof of the following theorem can be obtained by passing to the adjoint of the matrix representation of $\tilde{T}_A$ and using Theorem \ref{Cz}.

\begin{Theorem}
Let $A \colon\ \clm \raro \cln$ be a pure contraction. Then $\tilde{T}_A$ is a c.n.u. partial isometry and the characteristic function  of $\tilde{T}_A$ coincides with $\Theta (z) = A z^m$. Moreover, the upper triangular representation of $\tilde{T}_A$ in (\ref{tildeT_A}) coincides with the canonical representation introduced in Remark \ref{rem3.5}.
\end{Theorem}

\vspace{0.1in}

Concerning the operators $T_{0,b,c}$, it is easy to see that the representation (\ref{2.3}) is the canonical representation considered in Remark \ref{rem3.5} and therefore the quasi-similarity class in Theorem \ref{5.3} and Remark \ref{5.6} can be replaced with the dimension of $\clh_0$; thus for this very particular case, the answer to the problem in Remark \ref{5.6} is positive. The following theorem extends this observation from the scalar case to the vector valued case.

\begin{Theorem}\label{mono-N}
Let $T_A$ be the c.n.u. operator with monomial characteristic function $\Theta(z) = z^m A$. Then the nilpotent operators $N_0$ and $N_{*0}$ corresponding to the minimal representation considered in Theorem \ref{T5.2} are unitarily equivalent.
\end{Theorem}
\NI \textsf{Proof.} Let $\bm{h} = \cdots \oplus h_2 \oplus h_1 \oplus h_0 \oplus h_{-1} \oplus h_{-2} \oplus \cdots$ be a vector in the maximal $T_A$-invariant subspace $\clm_1$ of $\clh$ where $T_A$ is an isometry. First, observe that for all $1 \leq k < m$

\begin{equation*}
\begin{split}
T^k_A \bm{h} = & \cdots \oplus (D_{A^*} h_{m-(k-1)} - A h_{-(k-1)})_{m+1} \oplus (D_{A^*} h_{m-k} - A h_{-k})_{m} \oplus (h_{m-(k+1)})_{m-1} \oplus \cdots \\ & (h_0)_k \oplus (0)_{k-1} \oplus \cdots \oplus (0)_{1} \oplus (0)_0 \oplus (h_{-(k+1)})_{-1} \oplus (h_{-(k+2)})_{-2} \oplus \cdots
\end{split}
\end{equation*}

and
\begin{equation*}
\begin{split}
T^{m}_A \bm{h} = & \cdots \oplus (h_m)_{2m+1} \oplus (D_{A^*} h_{m-1} - A h_{-1})_{2m} \oplus \cdots \oplus (D_{A^*} h_{0} - A h_{-m})_{m} \\ & \oplus  (0)_{m-1} \oplus \cdots \oplus (0)_0 \oplus (h_{-m-1})_{-1} \oplus (h_{-m-2})_{-2} \oplus \cdots.
\end{split}
\end{equation*}
Thus $D^2_{T_A} T^{m}_A \bm{h} = 0$ yields $D_A h_{-m-1} = 0$ and hence $h_{-m-1} = 0$. Similarly, for all $k >m$ we have that $D^2_{T_A} T^k_A \bm{h} = 0$ and hence $h_{-k-1} = 0$. Then \begin{equation}\label{m1h}\bm{h} = \cdots \oplus h_2 \oplus h_1 \oplus h_0 \oplus h_{-1} \oplus \cdots \oplus  h_{-m} \oplus 0 \oplus 0 \oplus \cdots.\end{equation}
By virtue of $D^2_{T_A} T^k_A \bm{h} = 0$ for each $0 \leq k <m$ we obtain that $$A^* h_{m-(k+1)} + D_A h_{-(k+1)} = 0.$$ In other words, the elements of $\clm_1$ are precisely those $\bm{h} = \cdots \oplus h_2 \oplus h_1 \oplus h_0 \oplus h_{-1} \oplus \cdots \oplus  h_{-m} \oplus 0 \oplus 0 \oplus \cdots$ such that $A^* h_{m-(k+1)} + D_A h_{-(k+1)} = 0$ for all $k = 0, 1, \ldots, m-1$.
Set $g = D_{A^*} h_{m-(k+1)} - A h_{-(k+1)}$, then $h_{m-(k+1)} = D_{A^*} g, \, h_{-(k+1)} = - A^*g.$
Since for any given $g \in \cln$ the vectors $h_{m-(k+1)} = D_{A^*} g$ and $h_{-(k+1)} = - A^* g$ satisfy the relation $A^* h_{m-(k+1)} + D_A h_{-(k+1)} = 0$, we can infer that
\begin{equation}
\begin{split}
\clm_1 = & \clh_1 \oplus \{(D_{A^*} g_{m-1})_{m-1} \oplus (D_{A^*} g_{m-2})_{m-2} \oplus \cdots \oplus (D_{A^*} g_0)_{0} \oplus  (-A^* g_{-1})_{-1} \oplus \\& (-A^* g_{-2})_{-2} \oplus \cdots \oplus (-A^* g_{-m})_{-m} \oplus 0 \oplus 0 \oplus \cdots \colon\ g_i \in \cln, i = 0, 1, \ldots, m-1\}.
\end{split}
\end{equation}
Consequently
\begin{equation*}
\begin{split}
\clm_1 \cap \clh_0 = & \{(D_{A^*} g_{m-1})_{m-1} \oplus \cdots \oplus (D_{A^*} g_0)_0 \colon\ g_i \in \mbox{ker}\, A^*, i = 0, 1, \ldots, m-1\}\\ = & \{g_{m-1}\oplus \cdots \oplus g_0 \colon\ g_i \in \mbox{ker}\, A^*, i = 0, 1, \ldots, m-1\}.
\end{split}
\end{equation*}

Therefore, the minimal space $\clh_{00}$ is
\begin{equation*}
\begin{split}
\clh_{00} = & \clh_0 \ominus (\clm_1 \cap \clh_0)\\
= & (\mbox{ran}\, A)^- \oplus (\mbox{ran}\, A)^- \oplus \cdots \oplus (\mbox{ran}\, A)^-.
\end{split}
\end{equation*}
Therefore the matrix representation of the nilpotent operator $N_0$ is given by
\begin{equation*}
N_{0} = \begin{bmatrix} 0&I_{(\mbox{ran}\, A)^-}&0&\cdots &0&0&0\\0&0&I_{(\mbox{ran}\, A)^-}& \cdots & 0&0&0\\ 0& 0&0& \cdots & 0&0&0\\ \vdots & \vdots &\vdots& \cdots& \vdots &\vdots &\vdots \\ 0&0&0& \cdots & 0&I_{(\mbox{ran}\, A)^-}&0\\0 & 0 & 0&\cdots& 0&0&I_{(\mbox{ran}\, A)^-}\\0&0&0&\cdots &0&0&0\end{bmatrix}.
\end{equation*}

\NI Similarly, the other minimal space can be obtained as
\begin{equation*}
\begin{split}
\clh_{0*0} = & \clh_{0*} \ominus (\clh_{0*} \cap \clm_{-1})\\ = & \clh_{0*} \ominus \{ \{(D_{A} h_{-1})_{m-1} \oplus (D_{A} h_{-2})_{m-2} \oplus \cdots \oplus (D_{A} h_{-m})_{0} \colon\ h_{-k} \in \mbox{ker} A, k = 1, 2, \ldots, m\}\\ = & (\mbox{ran}\, A^*)^- \oplus (\mbox{ran}\, A^*)^- \oplus \cdots \oplus (\mbox{ran}\, A^*)^-,
\end{split}
\end{equation*}
and
\begin{equation*}
N_{*0} = \begin{bmatrix} 0&I_{(\mbox{ran}\, A^*)^-}&0&\cdots &0&0&0\\0&0&I_{(\mbox{ran}\, A^*)^-}& \cdots & 0&0&0\\ 0& 0&0& \cdots & 0&0&0\\ \vdots & \vdots &\vdots& \cdots& \vdots &\vdots &\vdots \\ 0&0&0& \cdots & 0&I_{(\mbox{ran}\, A^*)^-}&0\\0 & 0 & 0&\cdots& 0&0&I_{(\mbox{ran}\, A^*)^-}\\0&0&0&\cdots &0&0&0\end{bmatrix}.
\end{equation*}
It only remains to show that $\mbox{dim} (\mbox{ran}\, A^*)^- = \mbox{dim} (\mbox{ran}\, A)^-$. Indeed, with respect to the decompositions $\clm = (\mbox{ran} \,A^*)^- \oplus \mbox{ker}\, A$ and $\cln = (\mbox{ran} \,A)^- \oplus \mbox{ker}\, A^*$, we have
\begin{equation*}
A = \begin{bmatrix}A_{11}&0\\0&0\end{bmatrix},
\end{equation*}
where $A_{11} = P_{(\mbox{ran} \,A)^-} A|_{(\mbox{ran} \,A^*)^-}$ and $\mbox{ker}\, A_{11} = \{0\}$ and $\mbox{ker}\, A^*_{11} = \{0\}$. This shows that $\mbox{dim} (\mbox{ran} \,A^*)^- = \mbox{dim} (\mbox{ran}\, A)^-$, as required and hence, the Jordan forms $N_0$ and $N_{*0}$ are unitarily equivalent. The proof is complete. \qed

\vspace{0.2in}

We conclude with a simple observation on the spectrum of $T_A$ which follows from the relation of the spectrum $\sigma(T)$  of a contraction $T$ and the invertibility of the characteristic function $\Theta_T(z)$ (see \cite{NF}, page 259).

\begin{Remark}
The spectrum of the operator $T_A$ with the characteristic function  $\Theta(z)$ coinciding with $A z^m$ is $\partial \mathbb{D} \cup \{0\}$ if $A$ is invertible or $\overline{\mathbb{D}}$ if $A$ is not invertible.
\end{Remark}

We thank Professor J. A. Ball for kindly pointing to us that our work is related to realization theory for monic polynomials (see \cite{GLR} and \cite{R}). We note that the intersection of the class of polynomials considered in this paper with that of monic polynomials is the set of monomials $A z^n$ where $A = I$ for which the state operators, namely, $T_A^*$ are very particular.


\vspace{0.4in}

\begin{thebibliography}{99}



\bibitem{BM}
B. Bagchi and G. Misra, {\em Constant characteristic functions and homogeneous operators},  J. Operator Theory  37  (1997),  no. 1, 51--65. \MR98a:47009


\bibitem{B}
 H. Bercovici, {\em Operator theory and arithmetic in $H^{\infty}$}, Mathematical Surveys and Monographs, 26. American Mathematical Society, Providence, RI, 1988. \MR90e:47001

\bibitem{CM}
D. N. Clark and G. Misra, {\em On homogeneous contractions and unitary representations of SU$(1,1)$},  J. Operator Theory  30  (1993),  no. 1, 109--122. \MR96h:47015


\bibitem{CD}
M. J. Cowen and R. G. Douglas, {\em Complex geometry and operator theory}, Acta Math. 141 (1978), no. 3-4, 187--261. \MR{80f:47012}


\bibitem{FF}
C. Foias and A. Frazho, {\em The commutant lifting approach to interpolation problems}, Operator Theory: Advances and Applications, 44. Birkh\"auser Verlag, Basel, 1990  \MR92k:47033

\bibitem{GLR}
I. Gohberg, P. Lancaster and L. Rodman, {\em Matrix polynomials}, Academic Press, New York, 1982. \MR84c:15012

\bibitem{GGK}
I. C. Gohberg, S. Goldberg and M. Kaashoek, {\em Basic classes of linear operators}, Birkhauser Verlag, Basel, 2003. \MR2005g:47001

\bibitem{GK}
I. Gohberg and M. G. Krein, {\em Introduction to the theory of linear nonselfadjoint operators}, Translated from the Russian by A. Feinstein. Translations of Mathematical Monographs, Vol. 18 American Mathematical Society, Providence, R.I. 1969. \MR39 7447


\bibitem{NF}
B. Sz.-Nagy and C. Foias, {\em Harmonic Analysis of Operators on Hilbert Space}, North Holland, Amsterdam, 1970. \MR43 947.


\bibitem{Pa} S. Parrott, {\em On a quotient norm and the Sz.-Nagy -Foias lifting theorem}, J. Funct. Anal.  30  (1978), no. 3, 311--328. \MR81h:47006

\bibitem{Pe}
V. Peller, {\em Hankel operators and their applications}, Springer Monographs in Mathematics. Springer-Verlag, New York, 2003. \MR2004e:47040

\bibitem{Pi}
J. D. Pincus, {\em Commutators and systems of singular integral equations, I}, Acta Math. 121 (1968), 219-249. \MR39 2026

\bibitem{R}
L. Rodman, {\em An Introduction to Operator Polynomials}, OT 38, Birkhäuser Verlag, Basel and Boston, 1989. \MR90k:47032

\bibitem{T} R. Teodorescu, {\em Fonctions caract\'{e}ristiques constantes},
Acta Sci. Math. (Szeged) 38 (1976), no. 1-2, 183--185. \MR0442721

\bibitem{WP}
P.Y. Wu, {\em Contractions with constant characteristic function are reflexive},  J. London Math. Soc. (2)  29  (1984),  no. 3, 533--544. \MR0754939

\end{thebibliography}
\end{document}